\documentclass[a4paper,12pt]{article}
\RequirePackage[T1]{fontenc}
\RequirePackage{fix-cm}
\usepackage{amsmath,amsfonts,amsthm}
\usepackage{graphicx} 
\usepackage[unicode, final, colorlinks=true]{hyperref} 
\usepackage{cite}
\usepackage{caption}
\usepackage{mathtools}
\usepackage{newtxtext}
\usepackage[varvw]{newtxmath}
\usepackage{bm,bbm}
\usepackage[T1]{fontenc}
\usepackage[utf8]{inputenc}
\usepackage[dvipsnames]{xcolor}
\usepackage{tikz, pgfplots}
\usepackage{tikzscale}
\usepackage{booktabs}
\usepackage{multirow}
\usepackage{array}
\usepackage{csquotes}
\usepackage{indentfirst}
\usepackage{enumitem}
\usepackage{accents}

\hypersetup{citecolor=teal, linkcolor=BrickRed, urlcolor=NavyBlue}
\pgfplotsset{compat=1.17}
\usepgfplotslibrary{groupplots}
\usepgfplotslibrary{fillbetween}
\usetikzlibrary{backgrounds}
\usetikzlibrary{arrows.meta}
\pgfplotsset{plot coordinates/math parser=false}
\usetikzlibrary{external}
\usepgfplotslibrary{patchplots}
\usetikzlibrary{shapes.geometric}
\usetikzlibrary{backgrounds}
\usetikzlibrary{intersections}
\usetikzlibrary{spy}
\newtheorem{definition}{Definition}[section]
\newtheorem{theorem}{Theorem}[section]

\newtheorem{proposition}[theorem]{Proposition}
\newtheorem{remark}[theorem]{Remark}

\allowdisplaybreaks[4]

\newcommand{\R}{\mathbb{R}}
\newcommand{\N}{\mathbb{N}}
\newcommand{\Z}{\mathbb{Z}}

\newcommand{\defeq}{\coloneqq}

\newcommand{\eps}{\varepsilon}
\newcommand{\ddt}{\partial_t}
\newcommand{\dds}{\partial_s}
\newcommand{\ddx}{\partial_x}

\providecommand{\keywords}[1]{\textit{Keywords:} #1}
\providecommand{\msc}[1]{\textit{2020 MSC:} #1}

\title{Numerical schemes for coupled systems of nonconservative hyperbolic equations}
\author{
  Niklas~Kolbe$^{1}$\footnote{Corresponding author} \and
  Michael~Herty$^{1}$ \and 
  Siegfried~Müller$^1$
}

\date{
  \small
  $^1$Institute of Geometry and Applied Mathematics,\\ RWTH Aachen University, Templergraben 55,\\ 52062 Aachen, Germany\\
   \smallskip
   {\tt \{kolbe,herty,mueller\}@igpm.rwth-aachen.de} \\
   \smallskip
   \today
}
\begin{document}
\maketitle
\abstract{A new linear relaxation system for nonconservative hyperbolic systems is introduced, in which a nonlocal source term accounts for the nonconservative product of the original system. Using an asymptotic analysis the relaxation limit and its stability are investigated. It is shown that the path-conservative Lax--Friedrichs scheme arises from a discrete limit of an implicit-explicit scheme for the relaxation system. The relaxation approach is further employed to couple two nonconservative systems at a static interface. A coupling strategy motivated from conservative Kirchhoff conditions is introduced and a corresponding Riemann solver provided. A fully discrete scheme for coupled nonconservative products is derived and studied in terms of path-conservation. Numerical experiments applying the approach to a coupled model of vascular blood flow are presented.}


\msc{35L65, 35R02, 65M08}

\section{Introduction}
It is well-known, that the theory of weak and entropy solutions for systems of hyperbolic conservation laws, see e.g.~\cite{godlewski2021numer,dafermos2016hyper}, cannot be applied  to nonconservative problems in a straight-forward way. The main reason is that the product of a possibly discontinuous function with a measure is not well-defined. 
In the past decades several contributions, both theoretical and numerical have been proposed to tackle this problem. Without aiming to give here a complete review, we collect some results that will be used in the design of the numerical schemes below. In~\cite{dal1995defin,Ulbrich2003aa,BouchutJames1998aa,BouchutJames1999ad,LiuSandu2008aa} a concept for the treatment of nonconservative products is presented, which focuses on locally bounded and Borel-measurable functions with bounded variation. Therein, an interpretation of the nonconservative product as a real-valued signed Borel measure depending on a family of paths is presented, see also~\cite{volpert1967spaces,ambrosio2000funct}. By means of the Borel measures the concept of weak solutions and weak entropy solutions depending on the chosen family of paths is introduced, see~\cite{lefloch1988entrop,dal1995defin}. Applying the concept to  Riemann problems in~\cite{lefloch1989shock} it is verified that contact discontinuities corresponding to linearly degenerated fields do not depend on the path family. The influence of the path family on shocks is discussed in~\cite{abgrall2010}. Different approaches to deal with nonconservative products numerically are discussed in~\cite{castro2008why}. Those approaches have given rise to the class of path-conservative numerical methods including various finite volume~\cite{pares2006numer,castro2007well,castro2013hllc,castro2019path} and discontinuous Galerkin schemes~\cite{rhebergen2008discongaler}. 

We follow a similar concept in this paper, but are interested in the problem of {\em coupled } dynamics (in the spatially one-dimensional case). The coupling then induces a boundary condition at some static point in space, see Section~\ref{sec:4.1} for a precise definition. The problem of coupling spatially one-dimensional transport dynamics has gained attraction in the engineering and mathematical community in the past decades. It generalizes the related problem of conservation laws with discontinuous flux functions in space, see e.g.,~\cite{MR2396483, andreianov2011}, is core to flow problems on networks and enjoys a variety of applications, ranging from vehicular traffic on urban roads~\cite{GaravelloPiccoli2006ac}, gas and water flow in pipe systems~\cite{CorondAndrea-NovelBastin2007aa,BandaHertyKlar2006aa} to blood circulatory system in humans~\cite{Rathish-KumarQuateroniFormaggia2003aa}, as well as a wealth of mathematical models to represent and study other problems on directed graphs, see~\cite{BreCanGar2014} for a recent review and further references. In the mentioned references, {\em conservative} systems of hyperbolic balance laws have been studied leading to well-posedness results mostly based on Wave-Front-Tracking~\cite{Bressan:2000ab,Holden:2011aa}.
This construction relies heavily on the discussion of (half-) Riemann problems~\cite{HertyRascle2006aa,ColomboHertySachers2008aa} at the coupling interface. Therefore, extending those methods towards the coupling of nonconservative dynamics is not straight-forward. Coupled nonconservative dynamics appear in many applications, e.g., in the blood flow with dynamically changing cross-sections, see also Section~\ref{sec:blood} and references~\cite{MR2500548, formaggia2003one, quarteroni2004mathem}. Therein, the evolution of the cross section area is modeled by a nonconservative equation. Similar problems also appear in the modeling of multi-layer shallow water systems, see Section~\ref{sec:5.1} and \cite{castro2001}, or gas dynamics with variable pipe diameter~\cite{holle1}. In order to develop numerical schemes that may also be able to treat nonconservative dynamics, we follow a recent approach suggested in~\cite{herty2023centralsystems,herty2023centr,borsche2018kinet}. Therein, a scheme has been derived and analyzed that does {\em not} rely on the solution to nonlinear (half-) Riemann problems and does not require nonlinear Lax-curves. It employs a reformulation of the nonlinear system in a linear relaxation form proposed in~\cite{jin1995relaxschemsystem}. Due to the linearity of the relaxed systems the resulting coupling conditions could be explicitly solved (on a linear subspace computed by the constant eigenvectors of the system). Further, after taking the zero relaxation limit~\cite{jin1995relaxschemsystem}, an explicit numerical scheme for the coupled system is obtained.

In this paper, we follow a similar path to derive a numerical scheme, that eventually will be able to treat the coupling of nonconservative transport dynamics. In Section~\ref{sec:nonconservativesystems} we briefly review the theory of path-conservative schemes. This is required to discuss consistency of the proposed approach with existing results in the literature. The discussion therein is limited to the Cauchy problem. In Section~\ref{sec:relaxation}, we use a relaxation approach to embed the nonconservative system~\eqref{eq:nonconservativesystem} into a phase space model using stiff relaxation, see Equation~\eqref{eq:relsystem}. While relaxation has been considered in the context of nonconservative systems for specific applications~\cite{MR2491538,MR2588433} and to construct numerical schemes~\cite{MR2128209}, we propose a general relaxation system, in which the nonconservative product still appears, but only in the (stiff) source term. An asymptotic analysis shows the consistency with the formulation of Section~\ref{sec:nonconservativesystems}, see Theorem~\ref{thm:expansion}. Finally, we discretize the relaxation system in Section~\ref{sec:relaxationscheme} and drive the relaxation parameter in the discrete scheme to zero. This scheme can now be coupled at the interface; the details are given in Section~\ref{sec:coupling}. Based on this scheme, numerical experiments on the aforementioned problems are reported in Section~\ref{sec:num}. 

\section{Nonconservative systems and path-conservative schemes}\label{sec:nonconservativesystems}
In this section we set the preliminaries of our study by recalling the notions of nonconservative systems and path-conservative finite volume schemes. 
\subsection{Nonconservative hyperbolic systems and weak solutions}
We consider the first-order system
\begin{equation}\label{eq:nonconservativesystem}
\ddt U + A(U) \, \ddx U = 0 \quad \text{in}\quad \R^+ \times \R,
\end{equation}
with state variable $U(t,x)\in \Omega\subset \R^m$ for a smooth map $A:\Omega \rightarrow \R^{m \times m}$. We assume that the subset $\Omega$ is open and convex and that the system is strictly hyperbolic, i.e., the matrix $A(U)$ has $m$ real and distinct eigenvalues for all $U\in\Omega$. We suppose that initial data for the state in system~\eqref{eq:nonconservativesystem} at time $t=0$ is given by the function $U^0:\R \rightarrow \Omega$. 

We note that~\eqref{eq:nonconservativesystem} is a generalization of the system of conservation laws
\begin{equation}\label{eq:conservationsystem}
\ddt U + \ddx F(U) =0,  
\end{equation}
in which case $A(U)$ takes the role of the Jacobian of the smooth map $F:\Omega \rightarrow \Omega$, i.e., $A(U) =DF(U)$. A weak solution satisfies \eqref{eq:conservationsystem} in a distributional sense. In more details, let $U \in L^\infty(\R^+ \times \R)^m \cap BV(\R^+ \times \R)^m$ be piecewise regular with $D(t) \subset \R$ denoting the positions of its discontinuities. Then for $t>0$ the distribution $[\ddx F(U(t, \cdot))]$ is defined by
\begin{equation}\label{eq:conservativedistribution}
\begin{split}
  \langle [\ddx F(U(t, \cdot))], \phi \rangle &= \int_\R \phi(x)\,  \ddx F(U(t, x))\, dx \\
  &\quad + \sum_{x_D \in D(t)} \left(F(U(t, x_D^+)) - F(U(t, x_D^-))\right) \phi(x_D)
\end{split}
\end{equation}
for all test functions $\phi \in \mathcal{C}^\infty_0(\R)^m$ denoting by $x_D^-$ and $x_D^+$ the limit as $x$ approaches $x_D$ from the left and from the right, respectively. This distribution can be interpreted as a Borel measure. If $U$ is a weak solution of \eqref{eq:conservationsystem} then for all test functions $\phi \in \mathcal{C}^\infty_0(\R\times \R)^m$ it holds
\begin{equation}\label{eq:weaksolutionconservative}
- \int_{\R^+} \int_\R U(t, x) \, \ddt \phi(t,x) \, dx \, dt + \int_{\R^+}  \langle [\ddx F(U(t, \cdot))], \phi \rangle \, dt + \int_\R U^0(x) \, \phi(0, x) \, dx  = 0.
\end{equation}

In case of the nonconservative problem \eqref{eq:nonconservativesystem} the product $A(U)\, \ddx U$ is not well defined in a neighborhood around a disontinuity $x_D$. However, following \cite{dal1995defin, lefloch1999repres} an unambiguous definition can be given with respect to a selected family of paths.
\begin{definition}\label{def:path} A locally Lipschitz continuous map $\Phi:[0,1] \times \Omega \times \Omega \rightarrow \Omega$ is called a \emph{family of paths} iff
  \begin{enumerate}
  \item $\Phi(0, W^-, W^+) = W^-$ and $\Phi(1, W^-, W^+) = W^+$ for all $W^-$, $W^+ \in \Omega$,
  \item for any bounded set $\Omega_b\subset \Omega$ there are constants $C_1$ and $C_2$ such that it holds
    \begin{align*}
      \left| \dds \Phi(s; W^-, W^+)\right| &\leq C_1 |W^+-W^-| \\
      \left| \dds \Phi(s; W^-_1, W^+_1) - \dds \Phi(s; W^-_2, W^+_2)\right| &\leq C_2 \left( |W^-_1-W^-_2| +  |W^+_1-W^+_2| \right)
    \end{align*}
    for all $W^-, W^+, W^-_1, W^+_1, W^-_2,W^+_2 \in \Omega_b$ and almost every $s \in [0,1]$.
  \end{enumerate}
\end{definition}
Having selected the family of paths $\Phi$ we can interpret the nonconservative product $A(U)\, \ddx U$ by means of a Borel measure that we denote by $[A(U) \, \ddx U]_\Phi$. Considering again a piecewise regular solution $U$ this measure is defined in analogy to~\eqref{eq:conservativedistribution} as
\begin{equation}\label{eq:pathdistribution}
  \begin{split}
    \langle &[A(U(t, \cdot)) \, \ddx U(t, \cdot)]_\Phi, \phi \rangle  = \int_\R A(U(t, x))\,\ddx U(t, x) \phi(x) dx \\
      &+ \sum_{x_D \in D(t)} \left(\int_0^1 A(\Phi(s; U(t, x_D^-), U(t,x_D^+))) \dds \Phi(s; U(t, x_D^-), U(t,x_D^+)) \, ds \right) \phi(x_D)
  \end{split}
\end{equation}
for all test functions $\phi \in \mathcal{C}^\infty_0(\R)^m$. Replacing thus the distribution in~\eqref{eq:weaksolutionconservative} with~\eqref{eq:pathdistribution} gives rise to the concept of weak solution of~\eqref{eq:nonconservativesystem} with respect to the family $\Phi$.

\subsection{Path conservative schemes}\label{sec:pathconservativity}
In this section we describe generalized finite volume schemes for the nonconservative system~\eqref{eq:nonconservativesystem} and the path conservation property. To this end we fix the positive space and time increments $\Delta x$ and $\Delta t$ to define a uniform discretization of the time domain setting $t^n=n\Delta t$ and of the real line into the cells $I_j=[x_{j-1/2}, x_{j+1/2}]$, such that $x_{j-1/2}=j\,\Delta x$. Following the framework in \cite{pares2006numer} we obtain, in analogy to the conservative finite volume discretization, the path-conservative version
\begin{equation}\label{eq:pathconservationdiscretization}
  \begin{split}
  \frac{1}{\Delta x} \int_{I_j} U(x, t^{n+1}) \, dx &= \frac{1}{\Delta x} \int_{I_j} U(x, t^n) \, dx
    \\ &\quad - \frac{\Delta t}{\Delta x} \frac{1}{\Delta t} \int_{t^n}^{t^{n+1}}  \langle [A(U(t, \cdot)) \, \ddx U(t, \cdot)]_\Phi, \mathbbm{1}_{I_j} \rangle \, dt.
    \end{split}
\end{equation}
with $\mathbbm{1}$ denoting the indicator function. We introduce the cell averages $U_j^n$ in cell $I_j$ and at time $t_n$ and denote the piecewise constant approximation combining the cell averages $U_j^n$ at time $t^n$ over the full real line by $U^n$. Relying on an explicit discretization of the time integral in \eqref{eq:pathconservationdiscretization} we obtain
\begin{equation}\label{eq:notfullydiscretized}
  U_j^{n+1} = U_j^n - \frac{\Delta t}{\Delta x}  \langle [A(U^n)\, \ddx U^n]_\Phi, \mathbbm{1}_{I_j} \rangle
\end{equation}
as update formula. 
We note that in \eqref{eq:notfullydiscretized} the measure only consists of its singular constituent, i.e., weighted Dirac measures at the cell interfaces $x_{j-1/2}$ and $x_{j+1/2}$ due to $U^n$ being piecewise constant. Finite volume schemes can thus be formulated discretizing the product in~\eqref{eq:notfullydiscretized} by two terms corresponding to the contribution concerning both interfaces, i.e.,
\begin{equation}\label{eq:fvnonconservative}
    U_j^{n+1} = U_j^n - \frac{\Delta t}{\Delta x}  \left(D_{j-1/2}^{n,+} + D_{j+1/2}^{n,-}\right).
  \end{equation}
  These terms could, for example, be given as a function of the adjacent cell averages so that
\begin{equation}\label{eq:D}
    D_{j+1/2}^{n,\mp} = D^\mp(U_j, U_{j+1})
  \end{equation} given two maps $D^\mp:\Omega \times \Omega \to \R^n$.
  The scheme is called \emph{$\Phi$-conservative} if the map in \eqref{eq:D} is such that $D^\mp(U, U)=0$ for all $U\in \Omega$ and
  \begin{equation}\label{eq:phiconservative}
D^-(U_L, U_R) + D^+(U_L, U_R) =  \int_0^1 A(\Phi(s; U_L, U_R)) \dds \Phi(s; U_L, U_R) \, ds
\end{equation}
for all $U_L, U_R\in \Omega$.

An example of a path conservative scheme is the generalized form of the Lax-Friedrichs scheme, which is well known as numerical method for conservative hyperbolic systems of the form \eqref{eq:conservationsystem}. The generalized form is obtained by setting
\begin{equation}
D^\mp(U_L, U_R) = \frac 1 2 \int_0^1 A(\Phi(s; U_L, U_R)) \dds \Phi(s; U_L, U_R) \, ds \pm \frac{\Delta x}{2 \Delta t} (U_R - U_L),
\end{equation}
which results in the scheme
\begin{equation}\label{eq:laxfriedrich}
  \begin{split}
    U_j^{n+1} &= U_j^n - \frac{\Delta t}{2 \Delta x}  \int_0^1 A(\Phi(s; U_{j-1}^n, U_j^n)) \dds \Phi(s; U_{j-1}^n, U_j^n) \, ds \\
              &\quad -  \frac{\Delta t}{2 \Delta x}  \int_0^1 A(\Phi(s; U_{j}^n, U_{j+1}^n)) \dds \Phi(s; U_{j}^n, U_{j+1}^n) \, ds 
               + \frac 1 2 \left( U_{j-1}^n - 2 U_j^n + U_{j+1}^n\right) 
  \end{split}
\end{equation}
considered in~\cite{pares2009}.

\section{The relaxation system}\label{sec:relaxation}
Motivated by the approach in~\cite{jin1995relaxschemsystem} we consider a relaxation form of~\eqref{eq:nonconservativesystem}. To this end we select a family of paths $\Phi$ and introduce the auxiliary state $V(t, x)\in \R^m$ that together with the relaxation state $U$ is governed by the system
\begin{equation}\label{eq:relsystem}
  \begin{split}
    \ddt U + \ddx V &= 0\quad \text{in}\quad \R^+ \times \R, \\
    \ddt V + \Lambda \ddx U &= \frac{1}{\varepsilon} \left(  \langle [A(U(t, \cdot)) \, \ddx U(t, \cdot)]_\Phi, \mathbbm{1}_{(-\infty, x]} \rangle  - V \right)\quad \text{in}\quad \R^+ \times \R,
  \end{split}
\end{equation}
where $\varepsilon>0$ denotes the \emph{relaxation rate} and $\Lambda \in \R^{m \times m}$ is a diagonal matrix with positive entries $\lambda_1,\dots,\lambda_m$. Note that while both states of the system $U=U^\varepsilon$ and $V=V^\varepsilon$ depend on the relaxation rate, we neglect this dependency in the notation for simplicity. The variable $V$ is governed by a balance law with non-local source term employing the product~\eqref{eq:pathdistribution}. We consider the Cauchy problem with initial data given by $U^0$ and $V^0$ such that $V^0(x)=\langle [A(U^0) \, \ddx U^0]_\Phi, \mathbbm{1}_{(-\infty, x]} \rangle$.

\subsection{Asymptotic analysis}
To study the behavior of the solution of \eqref{eq:relsystem} as $\varepsilon \to 0$ we consider a Chapman-Enskog expansion assuming a representation of the states given by
\begin{equation}\label{eq:expansion}
  \begin{split}
    U &= U_0 + \varepsilon U_1 + \varepsilon^2 U_2 + \mathcal{O}(\varepsilon^3),\\
    V &= V_0 + \varepsilon V_1 + \varepsilon^2 V_2 + \mathcal{O}(\varepsilon^3)
  \end{split}
\end{equation}
such that for $j=0,1,2$ the functions $U_j$ and $V_j$ are piecewise smooth and, as well as their first spatial and temporal derivatives, independent of the relaxation rate $\varepsilon$. In analogy to Section~\ref{sec:nonconservativesystems} we collect the positions of discontinuities of all functions occurring in \eqref{eq:expansion} at time $t$ in the finite set $D(t)$. The statements given in the following analysis are to be understood for almost every $x\in \R$, i.e. for all $x\notin D(t)$.
\paragraph{Asymptotic analysis of the relaxation system} Inserting the expansion \eqref{eq:expansion} into the relaxation system \eqref{eq:relsystem} yields
\begin{equation}\label{eq:rel1sub}
    \ddt U_0 + \ddx V_0  + \varepsilon \left( \ddt U_1 + \ddx V_1 \right) = \mathcal{O}(\varepsilon^2) 
\end{equation}
for the first equation. To discuss the asymptotic expansion of the second equation we introduce the notation
\begin{equation}
T[U(t, \cdot)](x) \coloneqq \langle [A(U(t, \cdot)) \, \ddx U(t, \cdot)]_\Phi, \mathbbm{1}_{(-\infty, x]} \rangle
\end{equation}
and note that its $k$-th component is given by
\begin{equation}\label{eq:Tk}
  \begin{split}
    T^k&[U(t, \cdot)](x) = \int_{-\infty}^x \sum_{j=1}^m a_{kj}(U(t,y)) \ddx U^j \, dx  \\
    &\quad + \sum_{\substack{x_D \in D(t) \\ x_D \leq x}} \int_0^1 \sum_{j=1}^m a_{kj}(\Phi(s; U(t, x_D^-), U(t, x_D^+))) \dds \Phi^j(s; U(t, x_D^-), U(t, x_D^+)) \, ds
  \end{split}
\end{equation}
with upper indices indicating the component. Employing \eqref{eq:expansion} we expand this expression formally as
\begin{equation}
T^k[U] = T^k[U_0] + \varepsilon DT^k[U_0][U_1] + \mathcal{O}(\varepsilon^2),
\end{equation}
where $DT^k[U_0][U_1]$ denotes the Fréchet differential of $T^k$ at $U_0$ in direction $U_1$. An explicit form of this differential is considered below.  Let $DT[U_0][U_1]$ denote the vector with component $k$ given by $DT^k[U_0][U_1]$.  
Thanks to \eqref{eq:Tk} substituting~\eqref{eq:expansion} into the second equation of \eqref{eq:relsystem} yields
\begin{equation}\label{eq:rel2sub}
  \begin{split}
    \ddt V_0 &+ \Lambda \ddx U_0 
                 =  \frac{1}{\varepsilon} \left( T[U_0] - V_0 \right)
                 +  DT[U_0][U_1] - V_1+ \mathcal{O}(\varepsilon).
  \end{split}
\end{equation}
Comparing powers of $\varepsilon$ in both \eqref{eq:rel1sub} and \eqref{eq:rel2sub} gives rise to the two systems
\begin{equation}\label{eq:0order}
  \begin{split}
    \ddt U_0 + \ddx V_0 &= 0,\\
    V_0 &= T[U_0]
  \end{split}
\end{equation}
and
\begin{equation}\label{eq:1order}
\begin{split}
  \ddt U_1 + \ddx V_1 &= 0,\\
  \ddt V_0 + \Lambda \ddx U_0 &=  DT[U_0][U_1] - V_1 
\end{split}
\end{equation}
in terms of the state expansions of order zero and one. As the sum over the discontinuities in \eqref{eq:Tk} is independent of the space variable $x$ it holds
\begin{equation}\label{eq:ddxT}
\ddx T[U_0](t, x)  = A(U_0(t, x)) \ddx U_0(t, x).
\end{equation}
Thus, due to \eqref{eq:0order} the relaxation limit $U_0$ is a solution of the nonconservative system \eqref{eq:nonconservativesystem}, i.e.
\begin{equation}\label{eq:u0system}
\ddt U_0 + A(U_0)\ddx U_0 = 0.
\end{equation}

\paragraph{Asymptotic expansion of the relaxation system} Let $U$ and $V$ denote the solution of the relaxation system \eqref{eq:relsystem} and assume that it can be expanded as in \eqref{eq:expansion}. We note that component $k$ of the product $A(U) \ddx U$ can be expanded as
\begin{equation}\label{eq:compexpansion}
  \begin{split}
    \sum_{j=1}^m &a_{kj}(U) \ddx U^j = \sum_{j=1}^m \left(a_{kj}(U_0) +  \sum_{\ell=1}^m \varepsilon \partial_\ell a_{kj}(U_0) U_1^\ell \right) \ddx (U_0^j + \varepsilon U_1^j) + \mathcal{O}(\varepsilon^2)\\
    &=  \sum_{j=1}^m \left( a_{kj}(U_0) \ddx U_0^j + \varepsilon \left[ a_{kj}(U_0) \ddx U_1^j +  \sum_{\ell=1}^m \partial_\ell a_{kj}(U_0) U_1^\ell \ddx U_0^j \right] \right) + \mathcal{O}(\varepsilon^2)
  \end{split}
\end{equation}
with $\partial_\ell a_{kj}$ denoting the partial derivative with respect to component $\ell$ of the matrix-valued function $A$ restricted to its entry at row $k$ and column $j$. Inserting the expansion \eqref{eq:expansion} into \eqref{eq:nonconservativesystem} then yields
\begin{equation}\label{eq:nonconservativeexpansion}
  \begin{split}
    \ddt U + A(U) \, \ddx U &= \ddt U_0 + A(U_0) \, \ddx U_0  \\
    &\quad + \varepsilon \left( \ddt U_1 + A(U_0) \, \ddx U_1 + DA(U_0)(U_1) \, \ddx U_0 \right) + \mathcal{O}(\varepsilon^2),
  \end{split}
\end{equation}
where $DA(U_0)$ denotes the total differential of $A$ at $U_0$ mapping from $\Omega \subset \R^m$ to $\R^{m \times m}$. Applying this differential to $U_1$ and taking the product with $\ddx U_0$ a vector is obtained with $k$-th component given by
\[
(DA(U_0)(U_1) \, \ddx U_0)^k = \sum_{j=1}^m \sum_{\ell=1}^m \partial_\ell a_{kj}(U_0) U_1^\ell \ddx U_0^j.
  \]
 Employing the expansion~\eqref{eq:compexpansion} in~\eqref{eq:Tk} we note that
\begin{equation}
DT^k[U_0][U_1] = \int_{-\infty}^x  \sum_{j=1}^m a_{kj}(U_0) \ddx U_1^j +  \sum_{\ell=1}^m \partial_\ell a_{kj}(U_0) U_1^\ell \ddx U_0^j  \, dy + Q^k[U_0, U_1, \Phi]
\end{equation}
for a remainder term $Q^k$ that accounts for the discontinuities in~\eqref{eq:Tk} and satisfies $\ddx Q^k[U_0$, $U_1$, $\Phi]=0$ as it is independent of the spatial variable $x$. Consequently $\ddx DT[U_0][U_1]=A(U_0) \, \ddx U_1 + DA(U_0)(U_1) \, \ddx U_0$ and therefore due to the second equation of system~\eqref{eq:1order} it holds
\begin{equation}\label{eq:Atermidentity}
  A(U_0) \, \ddx U_1 + DA(U_0)(U_1) \, \ddx U_0 = \ddx \left(V_1 + \ddt V_0 +  \Lambda \ddx U_0 \right).
\end{equation}
Using the second equation of system~\eqref{eq:0order} we obtain
\begin{equation}
  \begin{split}
    \ddt V_0^k &= \int_{-\infty}^x \sum_{j=1}^m \left( a_{kj}(U_0) \ddt \ddx U^j_0 +  \ddx U^j_0 \sum_{\ell=1}^m \partial_\ell a_{kj}(U_0) \ddt U^\ell_0 \right) \, dx  + R^k[U_0, \Phi]
  \end{split}
\end{equation}
for component $k$, where the remainder term accounts for the temporal derivative of the term in~\eqref{eq:Tk} accounting for the discontinuities in $U_0$. This remainder term is again independent of the spatial variable $x$ and thus satisfies $\ddx R^k[U_0, \Phi]=0$.
Noting that due to~\eqref{eq:u0system} it holds $\ddt U_0 = -A(U_0) \ddx U_0$ we deduce using both the product and the chain rule 
\begin{equation}\label{eq:ddxddtV0}
  \begin{split}
    \ddx \ddt V_0 &= A(U_0) \ddt \ddx U_0 + DA(U_0)(\ddt U_0) \ddx U_0 \\
                  &= - A(U_0) \ddx(A(U_0) \ddx U_0) - DA(U_0)(A(U_0) \ddx U_0) \ddx U_0\\
                  &= - A(U_0) \ddx(A(U_0) \ddx U_0) - \ddx (A(U_0))A(U_0)\ddx U_0 \\
                  &\quad +  \ddx (A(U_0))A(U_0)\ddx U_0 -   DA(U_0)(A(U_0) \ddx U_0) \ddx U_0 \\
                  &=- \ddx \left( A(U_0)^2\ddx(U_0) \right) + DA(U_0)(\ddx U_0)(A(U_0) \ddx U_0)  \\
    &\quad -   DA(U_0)(A(U_0) \ddx U_0) \ddx U_0.
  \end{split}
\end{equation}
Simple algebra reveals that the latter difference can be written using a linear operator $M(U_0):\R^m\to \R^{m \times m}$ that is independent of the derivatives of $U_0$ such that
\begin{equation}\label{eq:M}
  M(U_0)(\ddx U_0) \ddx U_0 = DA(U_0)(\ddx U_0)(A(U_0) \ddx U_0)  -   DA(U_0)(A(U_0) \ddx U_0) \ddx U_0
\end{equation}
with components given by
\begin{equation}\label{eq:Mk}
  (M(U_0)(\ddx U_0) \ddx U_0)^k = \sum_{i=1}^m \sum_{j=1}^m \sum_{\ell=1}^m \ddx U_0^i \ddx U_0^j a_{\ell i}(U_0) \left(\partial_j a_{k\ell}(U_0) - \partial_\ell a_{kj}(U_0) \right).
\end{equation}

\noindent Substituting now \eqref{eq:Atermidentity} and \eqref{eq:ddxddtV0} in \eqref{eq:nonconservativeexpansion} we conclude
\begin{equation}
\begin{split}
  \ddt U + A(U) \ddx U &= \ddt U_0 + A(U_0) \ddx U_0  + \varepsilon \left( \ddt U_1 + \ddx V_1 \right) \\
  &\quad + \varepsilon \left( \Lambda \ddx^2 U_0 + \ddx \ddt V_0  \right) + \mathcal{O}(\varepsilon^2)\\
  &= \varepsilon \ddx \left( (\Lambda  - A(U_0)^2) \ddx U_0 \right) + \varepsilon M(U_0)(\ddx U_0) \ddx U_0 + \mathcal{O}(\varepsilon^2),
\end{split}
\end{equation}
where we have used \eqref{eq:u0system} and the first equation in system \eqref{eq:1order}. Using a suitable expansion of the linear operator $M(U)$ we infer
\begin{equation}
  \begin{split}
    \ddx \left( (\Lambda  - A(U)^2) \ddx U \right) + M(U)(\ddx U) \ddx U &=  \ddx \left( (\Lambda  - A(U_0)^2) \ddx U_0 \right)  \\
    &\quad + M(U_0)(\ddx U_0) \ddx U_0 + \mathcal{O}(\varepsilon),
  \end{split}
\end{equation}
which gives rise to the first order expansion that is summarized in the following theorem. 

\begin{theorem}[First order expansion]\label{thm:expansion}
 Suppose that system~\eqref{eq:nonconservativesystem} is hyperbolic and for a fixed family of paths $\Phi$ let $(U, V) \in L^\infty(\R^+ \times \R)^{2m} \cap BV(\R^+ \times \R)^{2m}$ be a weak solution of \eqref{eq:relsystem} that can be asymptotically expanded as in \eqref{eq:expansion}. Then, up to terms of order $\varepsilon^2$, the relaxation state $U$ is a weak solution of the equation
  \begin{equation}\label{eq:firstorderexpansion}
    \ddt U + A(U) \ddx U = \varepsilon \ddx \left( (\Lambda  - A(U)^2) \ddx U \right) + \varepsilon M(U)(\ddx U) \ddx U,
  \end{equation}
  where $M(U)(\ddx U) \ddx U$ is component-wise given by \eqref{eq:Mk}. In particular, as $\varepsilon \to 0$ the relaxation state $U$ solves the nonconservative system \eqref{eq:nonconservativesystem} in the weak sense.
\end{theorem}

The first order expansion \eqref{eq:firstorderexpansion} allows for a conclusion about the stability of the relaxation system.
\begin{remark}[Subcharacteristic condition]\label{rem:subchar}
 We assume that the elliptic terms are dominant and thus neglect the role of $M(U)(\ddx U) \ddx U$, which only includes first order derivatives of the state $U$, in the stability analysis of \eqref{eq:firstorderexpansion}. Then it is, however, necessary that
  \begin{equation}\label{eq:subcharacteristic}
    \Lambda - A(U)^2 \geq 0,
  \end{equation}
  i.e., the matrix $ \Lambda - A(U)^2$ should be positive semi-definite, for all $U \in \Omega$ to prevent anti-dissipative behavior of the relaxation system \eqref{eq:relsystem}. This requirement is consistent with the \emph{subcharacteristic condition} proposed in \cite{liu1987hyperconserlawsrelax} for the relaxation system from \cite{jin1995relaxschemsystem} in case of the conservative system~\eqref{eq:conservationsystem}. 
\end{remark}
In practice, the following proposition helps to find a suitable matrix $\Lambda$.
\begin{proposition}\label{prop:mu}
  Suppose that system \eqref{eq:nonconservativesystem} is hyperbolic and denote by
  \[
    \lambda_1(U) < \lambda_2(U) < \dots < \lambda_m(U)
  \]
  the eigenvalues of $A(U)$. Let further $\mu$ refer to the squared maximal spectral radius, i.e.,
  \[
    \mu = \max\{ \lambda_1(U)^2, \lambda_m(U)^2:~U\in \Omega \}.
  \]
  Then $\Lambda = \mu I$ satisfies the subcharacteristic condition \eqref{eq:subcharacteristic}.
\end{proposition}
\begin{proof}
Clearly, an eigenvector $v_k$ of $A(U)$ corresponding to the eigenvalue $\lambda_k(U)$ is an eigenvalue of $\mu I - A(U)^2$ corresponding to the eigenvalue $\mu - \lambda_k(U)^2$. Thus by the definition of $\mu$ the matrix $\Lambda - A(U)^2$ has $m$ distinct nonnegative eigenvalues for all $U\in \Omega$.
\end{proof}
We close the asymptotic analysis with a remark on the term~\eqref{eq:M}.
\begin{remark}
If problem \eqref{eq:nonconservativesystem} can be written in conservative form \eqref{eq:conservationsystem} it holds $M(U)=0$.
\end{remark}
\begin{proof}
In the conservative case it holds $A=DF$ for a smooth vector-valued flux function $F$, thus we have $\partial_j a_{k\ell} = \partial_j \partial_{\ell}F^k = \partial_\ell \partial_{j}F^k = \partial_\ell a_{kj}$ and the statements follows from the formulation \eqref{eq:Mk}.
\end{proof}

\subsection{The relaxation scheme}\label{sec:relaxationscheme}
Combining the first order upwind discretization with an implicit-explicit time discretization a scheme for the relaxation system~\eqref{eq:relsystem} is obtained. As \eqref{eq:relsystem} is a conservative system with nonlinear source term we can proceed in analogy to~\cite{hu2017asymppreserschem}, where a relaxation system for the conservative case was studied. As in Section \ref{sec:pathconservativity} we consider a discretization of the real line into uniform mesh cells $I_j=(x_{j-1/2}, x_{j+1/2})$ of width $\Delta x$ with origin located at the cell interface $x_{-1/2}$. Furthermore, the time line is partitioned into the instances $t^n=\sum_{k=1}^n \Delta t^n$ for some time increments $\Delta t^n>0$, which for brevity we assume uniform and denote by $\Delta t$ throughout this section. Let $U^n$ and $V^n$ denote a piecewise constant numerical solution of system~\eqref{eq:relsystem} in terms of cell averages at time $t^n$ and $U_j^n$ and $V_j^n$ the corresponding averages in cell $I_j$. The scheme admits the conservative form
\begin{equation}\label{eq:relschemeconservative}
  \begin{split}
  U^{n+1}_{j} &= U^{n}_{j} - \frac{\Delta t}{\Delta x}\left( F_{j+1/2}^{n} - F_{j-1/2}^{n} \right),\\
    V^{n+1}_{j} &= V^{n}_{j} - \frac{\Delta t}{\Delta x}\left( G_{j+1/2}^{n} - G_{j-1/2}^{n} \right)  + \frac {\Delta t} \eps \left( T_j[U^{n+1}]  - V^{n+1}_{j} \right)
  \end{split}
\end{equation}
for $j\in \Z$ with the numerical fluxes
\begin{equation} \label{eq:relflux}
  \begin{split}
  F_{j-1/2}^{n} =  \frac 1 2 \, (V_{j-1}^n + V_{j}^n)  - \frac 1 2  \sqrt{\Lambda} (U_j^n-U_{j-1}^n),\\
    G_{j-1/2}^{n} =  \frac 1 2 \Lambda  (U_{j-1}^n + U_{j}^n)  - \frac 1 2  \sqrt{\Lambda} (V_j^n-V_{j-1}^n).
  \end{split}
\end{equation}
By $\sqrt{\Lambda}$ we denote the diagonal matrix with positive entries $\sqrt{\lambda_i}$ for $i=1,\dots,m$. The term $T_j[U^n]$ discretizes the integral in~\eqref{eq:relsystem} and to account for discontinuities in the approximate solution makes use of the distributional form~\eqref{eq:pathdistribution} relying on a family of paths, i.e.,
  \begin{equation}\label{eq:T}
    T_j[U^{n}] \coloneqq \langle [A(U^n) \, \ddx U^n]_\Phi, \mathbbm{1}_{(-\infty, x_j]} \rangle
    = \sum_{i\leq j} \int_0^1 A(\Phi(s; U_{i-1}^n, U_i^n)) \dds \Phi(s; U_{i-1}^n, U_i^n)\, ds.
  \end{equation}
  In \eqref{eq:T} we use the notation $x_j$ for the center of the cell $I_j$, $\mathbbm{1}$ for the indicator function and we note that only the singular constituents of the distribution at the cell interfaces are considered as the approximate solution $U^n$ is constant in between.

\subsection{The relaxed scheme}\label{sec:relaxedscheme}
In this section, using asymptotic analysis we consider the relaxation limit of scheme~\eqref{eq:relschemeconservative}. The cell averages $U_j^n$ and $V_j^n$ occuring in the schemes above depend on the relaxation rate $\varepsilon$. We assume that for all $j$ and $n$ they can be asymptotically expanded around the zero-relaxation states $U_j^{n,0}$ and $V_j^{n,0}$ for sufficiently small relaxation rates as
\begin{equation}\label{eq:expansiondiscrete}
\begin{split}
  U_j^{n} = U_j^{n,0} + \varepsilon U_j^{n,1} + \mathcal{O}(\varepsilon^2),\\
  V_j^{n} = V_j^{n,0} + \varepsilon V_j^{n,1} + \mathcal{O}(\varepsilon^2).
\end{split}
\end{equation}
Combining the cell averages in the approximate solution $U^n$ at time $t^n$ we get the analogous expansion $U^n = U^{n,0} + \varepsilon U^{n,1} + \mathcal{O}(\varepsilon^2)$.
We fix $j\in \Z$ and then substitute the expansions~\eqref{eq:expansiondiscrete} into the second equation of~\eqref{eq:relschemeconservative}. After Taylor expansion we obtain
\begin{equation}\label{eq:vexpansion}
  \begin{split}
    V^{n+1, 0}_{j} \left(1 + \frac{\Delta t}{\varepsilon} \right) &= V^{n, 0}_{j} - \frac{\Delta t}{\Delta x}\left( G_{j+1/2}^{n, 0} - G_{j-1/2}^{n, 0} \right) 
     + \frac{\Delta t}{\varepsilon} \left( T_j[U^{n+1,0}]  \right)\\
    &\quad + \Delta t DT_j[U^{n+1,0}][U^{n+1,1}] - \Delta t V_j^{n+1,1}+ \mathcal{O}(\varepsilon),
    \end{split}
\end{equation}
where the term $G_{j-1/2}^{n, 0}$= $G_{j-1/2}^{n, 0}\left((U^{n,0}_{j-1}, V^{n,0}_{j-1}), (U^{n,0}_{j}, V^{n,0}_{j})  \right)$ denotes the numerical flux equivalent to $G_{j-1/2}$ in~\eqref{eq:relflux} at the zero-relaxation state
and $DT_j[U^{n+1,0}]:(\ell^\infty)^m \rightarrow \Omega$ is the Fréchet derivative of the operator~\eqref{eq:T} in $U^{n+1,0} \in (\ell^\infty)^m$.
We note that the following expansion holds
\begin{equation}
 \left( 1 + \frac{\Delta t}{\varepsilon} \right)^{-1} = \frac{\varepsilon}{\varepsilon + \Delta t} = \frac{\varepsilon}{\Delta t} + \mathcal{O}(\varepsilon^2), 
\end{equation}
and thus  by multiplication in~\eqref{eq:vexpansion} we obtain
\begin{equation}\label{eq:vlimitpre}
   V^{n+1, 0}_{j} =  T_j[U^{n+1,0}]   + \mathcal{O}(\varepsilon)
 \end{equation}
 for all $n\in \N_0$. Substituting now the asymptotic expansion~\eqref{eq:expansiondiscrete} into the first equation of~\eqref{eq:relschemeconservative} and taking into account~\eqref{eq:vlimitpre} we obtain
 \begin{equation}\label{eq:asymptoticlimitscheme}
   U^{n+1,0}_{j} = U^{n,0}_{j} - \frac{\Delta t}{\Delta x}\left( F_{j+1/2}^{n,0} - F_{j-1/2}^{n,0} \right)  + \mathcal{O}(\varepsilon),
 \end{equation}
 where the occurring numerical fluxes take the form
 \begin{equation}\label{eq:asymptoticfluxes}
   \begin{split}
     F_{j-1/2}^{n,0} =  \frac 1 2 \, ( V^{n, 0}_{j-1}+  V^{n, 0}_{j})  - \frac 1 2  \sqrt{\Lambda} (U_j^{n,0}-U_{j-1}^{n,0}).
   \end{split}
 \end{equation}
 In the relaxation limit the $\mathcal{O}(\varepsilon)$ terms in~\eqref{eq:vlimitpre} and ~\eqref{eq:asymptoticlimitscheme} vanish, which gives rise to the relaxed scheme
 \begin{equation}\label{eq:conservativeform}
  U_j^{n+1} = U_j^n - \frac{\Delta t }{\Delta x} \left( H_{j+1/2}^{n} -  H_{j-1/2}^{n}\right) \quad \text{for all }j \in \Z
\end{equation}
employing the numerical fluxes
\begin{equation}
  \begin{split}
  H_{j-1/2}^{n} &=  \frac 1 2 \, \left(T_{j-1}[U^{n}]  +  T_{j}[U^{n}]\right)  - \frac 1 2  \sqrt{\Lambda} (U_j^n-U_{j-1}^n).
  \end{split}
\end{equation}
In other words, the relaxed scheme reads
\begin{equation}\label{eq:relaxedscheme}
  \begin{split}
    U_j^{n+1} &= U_j^n 
 - \frac{\Delta t }{2\Delta x} \left( T_{j+1}[U^n] - T_{j-1}[U^n]\right)  + \frac{\Delta t }{2\Delta x}  \sqrt{ \Lambda} \left( U_{j-1}^n - 2 U_j^n + U_{j+1}^n\right),
  \end{split}
\end{equation}
where we note that
\begin{equation}
  \begin{split}
    T_{j+1}[U^n] - T_{j-1}[U^n] &= \int_0^1 A(\Phi(s; U_{j-1}^n, U_j^n)) \dds \Phi(s; U_{j-1}^n, U_j^n)\, ds\\
                                &\quad +  \int_0^1 A(\Phi(s; U_{j}^n, U_{j+1}^n)) \dds \Phi(s; U_{j}^n, U_{j+1}^n)\, ds.
  \end{split}
\end{equation}
Therefore, in the relaxation limit the nonlocality vanishes and an explicit local scheme is obtained.

\begin{proposition}
  The limit scheme \eqref{eq:relaxedscheme} is $\Phi$-conservative for any family of paths $\Phi$.
\end{proposition}
\begin{proof}
  Considering the terms \eqref{eq:D} and setting
  \[
    D^\mp(U_L, U_R) = \frac 1 2 \int_0^1 A(\Phi(s; U_L, U_R) \dds \Phi(s; U_L, U_R) \, ds \mp \frac{1}{2} \sqrt{\Lambda}(U_R - U_L)
  \]
  we can rewrite scheme \eqref{eq:relaxedscheme} in the form \eqref{eq:fvnonconservative}. Due to the second property in Definition~\ref{def:path} it holds $\dds \Phi(s; U, U)=0$ and therefore $D^\mp(U, U)=0$. Also property \eqref{eq:phiconservative} is easily verified.
\end{proof}

\begin{remark}
Assuming that the parameter matrix in the relaxation system \eqref{eq:relsystem} has the form $\Lambda = \lambda I$ and that $\Delta t$ is such that $\frac{\sqrt{\lambda} \Delta t}{\Delta x}=1$ the limit scheme \eqref{eq:relaxedscheme} coincides with the generalized Lax--Friedrichs scheme \eqref{eq:laxfriedrich}.
\end{remark}

\section{Coupling of nonconservative systems}\label{sec:coupling}
In this section we consider two nonconservative systems on the negative and the positive real half-axis, respectively:
\begin{equation}\label{eq:couplednonconservative}
  \begin{split}
    \ddt U + A_1(U) \, \ddx U &= 0 \quad \text{in}\quad \R^+ \times \R^-,\\
    \ddt U + A_2(U) \, \ddx U &= 0 \quad \text{in}\quad \R^+ \times \R^+.
  \end{split}
\end{equation}
We suppose that the state variable satisfies $U(t, x)\in \Omega_1$ for $x<0$ and $U(t, x)\in \Omega_2$ for $x>0$ for two open and convex subsets $\Omega_1 \subset \R^{m_1}$ and $\Omega_2 \subset \R^{m_2}$. In this setting the matrix-valued functions $A_i:\Omega_i \rightarrow \R^{m_i\times m_i}$ for $i\in\{1,2\}$ determine the system on the left and the right half-axis, respectively. We assume that restricted to either side the system is hyperbolic and that both $U-U_a^1$ is of compact support on the negative half-axis and $U-U_a^2$ is of compact support on the positive half-axis for fixed $U_a^i\in\Omega_i$ and $i\in\{1,2\}$. At the interface located at the origin the coupling condition
\begin{equation}\label{eq:psiu}
\Psi_U(U(t, 0^-), U(t, 0^+))=0 \qquad \text{for a. e. }t>0
\end{equation}
is imposed for a suitable mapping $\Psi_U:\R^{m_1} \times \R^{m_2} \rightarrow \R^p$. The number of conditions $p \in \N$ depends here on the coupled systems at hand. Initial data for system \eqref{eq:couplednonconservative} is given by the vector-valued function $U^0(x)$ taking values in $\Omega_1$ on the left and in $\Omega_2$ on the right half axis such that at the interface it holds $\Psi_U(U^0(0^-), U^0(0^+))=0$.

\subsection{Coupled relaxation system}\label{sec:4.1}

\begin{figure}
  \centering 
  \includegraphics{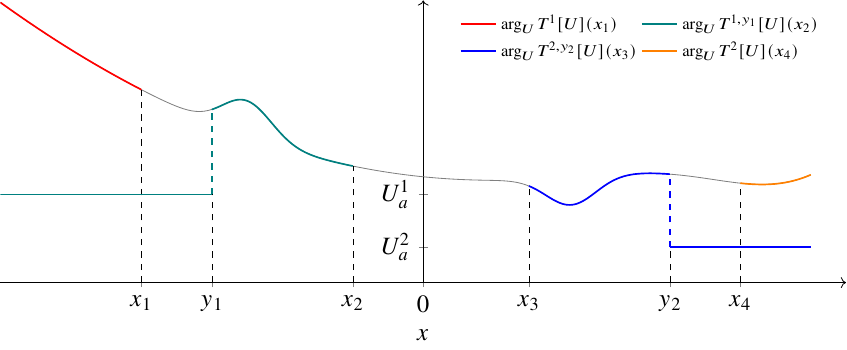}
    \caption{The function argument of $T^{i,y}[U](x)$ for $i\in\{1,2\}$ and various choices of $x$ and $y$ in an example case with a smooth scalar function $U$ (gray line) representing the solution of~\eqref{eq:coupledrelaxation} at a fixed time instance. Note that for finite $y$ a jump from $U_a^1$ or to $U_a^2$ is taken into account.}\label{fig:T}
    \end{figure}

In analogy to Section~\ref{sec:relaxation} we consider a relaxation system for the coupled problem~\eqref{eq:couplednonconservative}. To define suitable coupling conditions we use for any interval $I\subset \R$ the notation
\begin{equation}
  U_I(t,x) \coloneqq
  \begin{cases*}
    U(t,x) & if $x\in I$\\
    U_a^1 & if $x \notin I$ and $x<0$\\
    U_a^2 & if $x \notin I$ and $x>0$\\
  \end{cases*}                  
\end{equation}
to refer to the truncated state variable. Given the family of paths $\Phi$ we vary the direction of integration and define for any $y < x < 0$ the non-local operator
\begin{equation}\label{eq:T1}
  T^{1, y}[U(t, \cdot)](x) \coloneqq \langle [A_1(U_{(y, 0)}(t, \cdot) \, \ddx U_{(y, 0)}(t, \cdot)]_\Phi,  \mathbbm{1}_{(-\infty, x]} \rangle
\end{equation}
for the left half-axis and for any $y >x>0$ the operator
\begin{equation}\label{eq:T2}
T^{2,y}[U(t, \cdot)](x) \coloneqq - \langle [A_2(U_{(0, y)}(t, \cdot) \, \ddx U_{(0, y)}(t, \cdot)]_{\Phi}, \mathbbm{1}_{[x, \infty)} \rangle
\end{equation}
for the right half-axis as well as the notations $T^1\coloneqq T^{1,-\infty}$ and $T^2 \coloneqq T^{2, \infty}$.
For simplicity we often neglect the arguments in the operators \eqref{eq:T1} and \eqref{eq:T2}. Figure~\ref{fig:T} visualizes the function arguments of the operators \eqref{eq:T1} and \eqref{eq:T2} defined as $\arg_U T^{1,y}[U](x) \coloneqq U_{(y,0)}(\cdot) \rvert_{(-\infty, x]}$ and $\arg_U T^{2,y}[U](x) \coloneqq U_{(0,y)}(\cdot) \rvert_{[x, \infty)}$. We introduce the relaxation system
\begin{equation}\label{eq:coupledrelaxation}
  \begin{split}
    \ddt U + \ddx V &= 0\quad \text{in}\quad \R^+ \times \R \setminus \{ 0\}, \\
    \ddt V + \Lambda_1 \ddx U &= \frac{1}{\varepsilon} \left( T^1[U] - V \right)\quad \text{in}\quad \R^+ \times \R^-,\\
    \ddt V + \Lambda_2 \ddx U &= \frac{1}{\varepsilon} \left( T^2[U] - V \right)\quad \text{in}\quad \R^+ \times \R^+
  \end{split}
\end{equation}
governing the relaxation state $U$ and the auxiliary variable $V$ both depending on the relaxation rate $\varepsilon>0$. We introduce the vector $Q=(U,V)$ combining the two variables. By $\Lambda_1$ and $\Lambda_2$ we refer to diagonal matrices with positive entries that, as discussed in Proposition~\ref{rem:subchar}, satisfy subcharacteristic conditions for stability, i.e.,
\begin{equation}
\Lambda_i - A_i(U)^2 \geq 0 \quad \text{for all } U \in \Omega_i
\end{equation}
and $i=1,2$. At the interface we consider a mapping $\Psi_Q:\R^{2m_1}\times \R^{2m_2}\rightarrow \R^q$ for $q\in \N$ and impose the relaxation coupling condition
\begin{equation}\label{eq:relaxationcoupling}
  \Psi_Q(Q(t, 0^-), Q(t, 0^+)) = 0,
\end{equation}
which closes the coupled PDE problem. The initial data for the coupled relaxation system is adapted from the original initial data setting $Q^0=\left(U^0, V^0\right)$, so that
\begin{equation}
V^0(x)= T^1\left[U^0\right](x)\quad \text{if }x<0, \quad V^0(x)=T^2\left[U^0\right](x)\quad \text{if }x>0
\end{equation}
and compatibility of the initial data $Q^0$ with the coupling condition~\eqref{eq:relaxationcoupling} is given. Clearly, as we restrict system \eqref{eq:coupledrelaxation} to either the left or the right half-axis it satisfies the limit property, see Theorem~\ref{thm:expansion}. To study the relaxation limit at the interface we follow our framework in \cite{herty2023centralsystems} for conservative systems. An asymptotic expansion at the interface motivates the following notion of consistency.

\begin{definition}\label{def:consistency}
  The coupled relaxation system~\eqref{eq:coupledrelaxation} is \emph{consistent} with the nonconservative coupled system~\eqref{eq:couplednonconservative} iff for a.~e. $t>0$ the corresponding coupling conditions~\eqref{eq:relaxationcoupling} and~\eqref{eq:psiu} satisfy 
  \begin{equation}
    \Psi_U\left(U(t, 0^-), U(t, 0^+)\right) = 0 \quad \text{iff} \quad \Psi_Q \left(
    \begin{pmatrix} U(t, 0^-) \\
      T^{1}[U(t, \cdot)](0^-)
    \end{pmatrix},
    \begin{pmatrix} U(t, 0^+) \\
      T^{2}[U(t, \cdot)](0^+)
    \end{pmatrix}
  \right) = 0.
\end{equation}
\end{definition}

\begin{remark}\label{rmk:couplingconditions}
  Our approach allows for a generalized notion of the coupling conditions \eqref{eq:psiu} and \eqref{eq:psiq}: defining the two functionals  
  \begin{equation*}
    S^1: L^\infty(\R^-)^{m_1} \to C^0(\R^-)^{m_1}, \qquad  S^2: L^\infty(\R^+)^{m_2} \to C^0(\R^+)^{m_2}
  \end{equation*}
  we can consider the coupling condition
  \begin{equation}
    \Psi_U(S^1[U(t, \cdot)](0^-), S^2[U(t, \cdot)](0^+)) =0.
  \end{equation}
  In this setting, $S^1$ and $S^2$ depend only on $U$ restricted to the left and right half-axis, respectively. We can analogously consider coupling conditions of the type
  \begin{equation}
    \Psi_Q(\tilde S^1[Q(t, \cdot)](0^-), \tilde S^2[Q(t, \cdot)](0^+)) =0,
  \end{equation}
  where in the definition of the functionals $\tilde S^1$ and $\tilde S^2$ the number of components is doubled compared to $S^1$ and $S^2$. In fact, the case $S^1=T^{1, 0^-}$, $S^2=T^{2, 0^+}$ is discussed below.
\end{remark}
\begin{remark}\label{rmk:operatorrelation}
  The operator \eqref{eq:T1} satisfies for any $x<0$ the property
  \begin{multline*}
    T^1[U](y) +  T^{1,y}[U](x) - T^1[U](x) \\
    = \int_0^1 A_1(\Phi(s; U_a^1, U(y))) \dds \Phi(s; U_a^1, U(y)) \, ds \quad \text{for a.e. }y<x
  \end{multline*}
  and the operator \eqref{eq:T2} for any $x>0$ the property 
    \begin{multline*}
      T^{2,y}[U](x)
      + T^2[U](y)
      - T^2[U](x) \\
    = - \int_0^1 A_2(\Phi(s; U(y), U_a^2)) \dds \Phi(s;U(y),U_a^2) \, ds \quad \text{for a.e. }y>x.
  \end{multline*}
  Therefore, we obtain in the limit
  \begin{align*}
    \lim_{\delta \nearrow \, x}  \int_0^1 A_1(\Phi(s; U_a^1, U(\delta))) \dds \Phi(s; U_a^1, U(\delta)) \, ds  - T^1[U](\delta) +  T^1[U](x), \\
    \lim_{\delta \searrow \, x} - \int_0^1 A_2(\Phi(s; U(\delta), U_a^2)) \dds \Phi(s;U(\delta),U_a^2) \, ds  - T^2[U](\delta) + T^2[U](x),
  \end{align*}
  which both simplify to a single (signed) path integral if $U$ is continuous in $x$. In the following the latter expressions are referred to by $T^{1,x}[U](x)$ and $T^{2, x}[U](x)$, respectively.
\end{remark}

\subsection{(Half-) Riemann solvers for linear problems}
A Riemann solver (RS) for system \eqref{eq:coupledrelaxation} identifies suitable boundary/coupling data that solves the two half-Riemann problems at the interface, see~\cite{dubois1988boundconditnonlin}, and satisfies the relaxation coupling condition \eqref{eq:relaxationcoupling}. As~\eqref{eq:coupledrelaxation} is a {\em linear first order system}, it has a simple eigenvalue/eigenvector structure. Therefore, the solution does not require to compute (nonlinear) Lax-curves, but is directly given within a linear subspace. Those spaces are introduced below:
\begin{subequations}\label{eq:laxcurves}
\begin{equation}
 \mathcal{L}^-(Q_0^*) \coloneqq
  \left\{
    \begin{pmatrix}
      U_0^* -  (\sqrt{\Lambda_1})^{-1} \Sigma^- \\
      V_0^* + \Sigma^-
    \end{pmatrix} \in \R^{2m_1}\,:
    \Sigma^- \in \R^{m_1}
  \right\}
\end{equation}
comprises all states that connect to $Q_0^*=(U_0^*, V_0^*)$ by Lax-curves with negative speeds. Similarly, the space
\begin{equation}
  \mathcal{L}^+(Q_0^*) \coloneqq 
  \left\{
    \begin{pmatrix}
      U_0^* + (\sqrt{\Lambda_2})^{-1} \Sigma^+ \\
      V_0^* + \Sigma^+
    \end{pmatrix}\in \R^{2m_2}\,:
    \Sigma^+ \in \R^{m_2}
  \right\}
\end{equation}
\end{subequations}
contains all states connecting to $Q_0^*$ by Lax-curves with positive speeds. Note that both spaces are linear in the parametrization $\Sigma^\pm$, respectively. This is the major advantage also from a numerical point of view compared to the original nonconservative problem, where such a characterization does not exist.

Suppose that we are given discrete piecewise constant data next to the coupling interface from a numerical scheme. This so-called trace-data is denoted $Q_0^- = (U_0^-, V_0^-)$ left from the interface and $Q_0^+ = (U_0^+, V_0^+)$ right from the interface, respectively. The sought boundary/coupling data is referred to by the notation $Q_R=(U_R, V_R)$ in case of the left half-axis and by $Q_L=(U_L, V_L)$ in case of the right half-axis. The data $Q_R$ and $Q_L$ solve the half-Riemann problems at the interface iff the conditions
\begin{equation}\label{eq:laxcondition}
  Q_R \in \mathcal L^-(Q_0^-) \qquad \text{and} \qquad Q_L \in \mathcal L^+(Q_0^+)
\end{equation}
hold. A RS for the relaxation system~\eqref{eq:coupledrelaxation} and coupling condition~\eqref{eq:relaxationcoupling} is a mapping
\begin{equation}\label{eq:RSdef}
\mathcal{RS}: \R^{2m_1 \times 2m_2}\rightarrow \R^{2m_1 \times 2m_2}, \qquad (Q_0^-, Q_0^+) \mapsto (Q_R, Q_L)
\end{equation}
that assigns coupling data satisfying both \eqref{eq:laxcondition} and the coupling condition
\begin{equation}\label{eq:psiq}
\Psi_Q(Q_R, Q_L) = 0.
\end{equation}
We refer to~\cite{herty2023centralsystems} for a discussion of the well-posedness of those linear coupling problem given by \eqref{eq:laxcondition} and \eqref{eq:psiq} and note that in case of multiple solutions additional problem specific criteria are required to define a suitable RS.

Let $U(t, \cdot)$ be the first component of a piecewise constant solution to \eqref{eq:coupledrelaxation} at the time $t$ with trace data $U_0^-$ and $U_0^+$ left and right from the interface and assigned coupling data $U_R$ and $U_L$. Then the coupling states of $V$ in the relaxation limit required for consistency in Definition~\ref{def:consistency} take the form
\begin{align}
  T^{1}[U](0^-) &= T^{1}[U]( - \delta) + \int_0^1 A_1(\Phi(s; U_0^-, U_R)) \dds \Phi(s; U_0^-, U_R) \, ds,  \label{eq:T1interface} \\
  T^{2}[U](0^+) &= T^{2}[U]( \delta ) - \int_0^1 A_2(\Phi(s; U_L, U_0^+)) \dds  \Phi(s; U_L, U_0^+) \, ds  \label{eq:T2interface}
\end{align}
for small $\delta>0$ such that $U(t, -\delta) = U_0^-$ and $U(t, -\delta) = U_0^+$. The vectors $T^1[U]( - \delta)$ and $T^2[U](\delta)$ constitute the trace data for the auxiliary variable $V$ in the relaxation limit.

\subsubsection{Path conservative Kirchhoff conditions}\label{sec:kirchhoff}
As an example we discuss a RS for \emph{path conservative Kirchhoff conditions} in the case $m_1=m_2$, which we define such that the coupling function complementing system \eqref{eq:couplednonconservative} is given by
\begin{equation}\label{eq:kirchhoff}
\Psi_U =  T^{2, 0^+}[U](0^+) - T^{1, 0^-}[U](0^-).
\end{equation}
The occurring limits are taken with respect to the integration parameter $x$ in the operator introduced in Remark~\ref{rmk:operatorrelation}, i.e.,
\[
  T^{1, 0^-}[U](0^-) = \lim_{x \nearrow \, 0 } T^{1, x}[U](x), \qquad   T^{2, 0^+}[U](0^+) = \lim_{x \searrow \, 0} T^{2, x}[U](x).
  \]
  We consider the corresponding relaxation system \eqref{eq:coupledrelaxation} and employ Remark~\ref{rmk:operatorrelation} to derive the conditions
  \begin{align}
  T^{1, 0^-}[U](0^-) &= T^{2, 0^+}[U](0^+), \label{eq:relaxationkirchhoff1}\\
  \lim_{x \nearrow \, 0} ~\lim_{\delta \nearrow \, x }~ V(x) - V(\delta) + \int_0^1& A_1(\Phi(s; U_a^1, U(\delta))) \dds \Phi(s; U_a^1, U(\delta)) \, ds \notag \\
  = \lim_{x \searrow \, 0}~\lim_{\delta \searrow \, x} ~ V(x) - V(\delta) &- \int_0^1 A_2(\Phi(s; U(\delta), U_a^2)) \dds \Phi(s; U(\delta), U_a^2) \, ds \label{eq:relaxationkirchhoff2}
  \end{align}
  determining the coupling function $\Psi_Q$ so that the coupled relaxation system is consistent with the nonconservative coupled system.

  In the following we construct a RS for the relaxation system \eqref{eq:coupledrelaxation} with coupling conditions \eqref{eq:relaxationkirchhoff1} and \eqref{eq:relaxationkirchhoff2}. To this end we need to find the coupling data $U_R$, $U_L$, $V_R$, $V_L$ for given trace data $U_0^-$, $U_0^+$ and $V_0^-$ and $V_0^+$.
  The piecewise constant information next to the boundary allows for an interpretation of \eqref{eq:relaxationkirchhoff2} in terms of trace and boundary data motivating the definition of the quantities 
  \begin{align}
    \widetilde V_R &\coloneqq V_R - V_0^- + P_1(U_0^-), \qquad P_1(U_0^-) \coloneqq \int_0^1 A_1(\Phi(s; U_a^1, U_0^-)) \dds \Phi(s; U_a^1, U_0^-) \, ds, \label{eq:modVR}\\
    \widetilde V_L &\coloneqq V_L - V_0^+ - P_2(U_0^+), \qquad P_2(U_0^+) \coloneqq \int_0^1 A_2(\Phi(s; U_0^+, U_a^2)) \dds \Phi(s; U_0^+, U_a^2) \, ds, \label{eq:modVL}
  \end{align}
  of which $\widetilde V_R$ and $\widetilde V_L$ are discrete representations of the left- and right-hand side of \eqref{eq:relaxationkirchhoff1}. Employing the parametrization~\eqref{eq:laxcurves} as well as the identities~\eqref{eq:T1interface} and \eqref{eq:T2interface} we obtain for~\eqref{eq:relaxationkirchhoff1} and \eqref{eq:relaxationkirchhoff2} the system
  \begin{align}
    P_1(U_0^-) + \Sigma^- &= -  P_2(U_0^+) + \Sigma^+, \label{eq:rssystem1}\\
    P_1(U_0^-)  + \int_0^1 A_1(\Phi(s; U_0^-, U_0^- &-  (\sqrt{\Lambda_1})^{-1} \Sigma^-)) \dds \Phi(s; U_0^-, U_0^- -  (\sqrt{\Lambda_1})^{-1} \Sigma^-) \, ds \notag \\
    = -  P_2(U_0^+) - \int_0^1 A_2(\Phi(s; U_0^+ &+ (\sqrt{\Lambda_2})^{-1} \Sigma^+, U_0^+)) \dds  \Phi(s; U_0^+ + (\sqrt{\Lambda_2})^{-1} \Sigma^+, U_0^+) \, ds. \label{eq:rssystem2}
  \end{align}
  Any solution of this system for $(\Sigma^-, \Sigma^+)$ determines suitable coupling data after using the parametrization \eqref{eq:laxcondition} again. Provided, a solution to \eqref{eq:rssystem1} and \eqref{eq:rssystem2} exists a RS is defined specifying a suitable solution for the application at hand. 

  \begin{remark}
    Suppose that the two systems in \eqref{eq:couplednonconservative} are conservative so that $A_1=DF_1(U)$ and $A_2=DF_2(U)$. If $F_1(U_a^1)= F_2(U_a^2)$ it is easy to verify that given \eqref{eq:kirchhoff} the condition $\Psi_U =0$ is equivalent to the classical form of the Kirchhoff condition $F_1(U(0^-))= F_2(U(0^+))$. Otherwise the classical Kirchhoff condition is equivalent to the modified condition
    \begin{equation}\label{eq:modkirchhoff}
      T^{1, 0^-}[U](0^-) + F_1(U_a^1) =  T^{2, 0^+}[U](0^+) + F_2(U_a^2),
    \end{equation}
    for which the RS above can be adapted by adding $F_1(U_a^1)$ to the left-hand sides and $F^2(U_a^2)$ to the right-hand sides of \eqref{eq:rssystem1} and \eqref{eq:rssystem2}.
  \end{remark}
  
\begin{remark}
  While the coupling function
  \[\Psi_U^\prime=  T^{2}[U](0^+) - T^{1}[U](0^-)\]
  also leads to a coupling condition that under the assumption $F_1(U_a^1)= F_2(U_a^2)$ is equivalent to the classical Kirchhoff condition in the conservative case, this choice leads to a non-local RS for nonconservative systems causing artifacts even in case of a constant solution at the interface.
\end{remark}

\subsection{The relaxed scheme for coupled problems}
\begin{figure}
  \centering
  \includegraphics{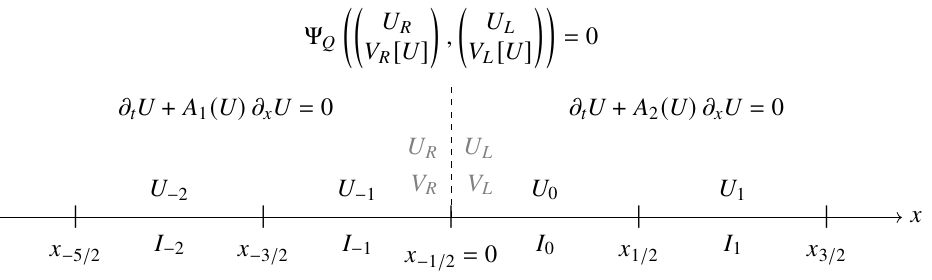}
\caption{The coupled non-conservative system as limit of the coupled relaxation system with consistent coupling condition on the discretized real line. The coupling data $V_R$ and $V_L$ are obtained from the numerical solution on the respective half-axis by~\eqref{eq:T1interface} and~\eqref{eq:T2interface}.}\label{fig:numcoupling}
\end{figure}
Combining a RS for the coupled relaxation system with the limit scheme derived in Section~\ref{sec:relaxedscheme} we derive a new fully-discrete scheme for the coupled nonconservative system~\eqref{eq:couplednonconservative}. We consider a relaxation system with coupling function $\Psi_Q$ consistent to the coupled problem and assume that a well-defined and continuously differentiable RS corresponding to $\Psi_Q$ is given. We adapt the discretization introduced in Section~\ref{sec:relaxationscheme} and impose the coupling interface between the mesh cells $I_{-1}$ and $I_0$ as indicated in Figure~\ref{fig:numcoupling}. Moreover, we introduce the discrete forms of the operators $T^1$ and $T^2$ for $j\in \Z$ as
\begin{align}
  T_j^1[U^{n}] &\coloneqq T^1[U^n](x_j)  = \sum_{i\leq j} \int_0^1 A_1(\Phi(s; U_{i-1}^n, U_i^n)) \dds \Phi(s; U_{i-1}^n, U_i^n)\, ds, \label{eq:T1discrete}
  \\    T_j^2[U^{n}] &\coloneqq T^2[U^n](x_j) = -\sum_{i\geq j} \int_0^1 A_2(\Phi(s; U_{i-1}^n, U_i^n)) \dds \Phi(s; U_{i-1}^n, U_i^n)\, ds. \label{eq:T2discrete}
\end{align}

Following the limit procedure in Section~\ref{sec:relaxedscheme} and making use of the smoothness of the RS the limit scheme for the coupled problem is derived. We also refer to \cite{herty2023centralsystems}, where the detailed computation is given for the conservative case. Adapting the common form \eqref{eq:fvnonconservative} we end up with the scheme
\begin{equation}\label{eq:relaxedcouplingscheme}
  U_j^{n+1} = U_j^n - \frac{\Delta t }{\Delta x} \left( D_{j+1/2}^{n,-} +  D_{j-1/2}^{n,+}\right) \quad \text{for all }j \in \Z,
\end{equation}
where for $j \in \Z \setminus \{ 0\}$, i.e., away from the coupling, the interface contributions are given by
\begin{align}
D_{j-1/2}^{n,\mp} = \frac{1}{2}\left( T_j^i[U^n] - T_{j-1}^{i}[U^n] \right)  \pm \frac{1}{2} \sqrt{\Lambda_i} \left( U_j^n - U_{j-1}^n \right)
\end{align}
with $i=1$ if $j<0$ and $i=2$ if $j>0$. At the interface we have
\begin{align}
  D_{-1/2}^{n,+} &=  \frac{1}{2}\left( V_R^n - T_{-1}^{1}[U^n] \right)  - \frac{1}{2} \sqrt{\Lambda_1} \left( U_R^n - U_{-1}^n \right),\\
  D_{-1/2}^{n,-} &=  \frac{1}{2}\left( T_0^2[U^n] - V_L^n \right)  + \frac{1}{2} \sqrt{\Lambda_2} \left( U_0^n - U_{L}^n \right) 
\end{align}
for coupling data given by 
\begin{equation}\label{eq:RSlimit}
  \mathcal{RS}( (U_{-1}^n, T_{-1}^1(U^n)), (U_0^n, T^2_0(U^n))) \eqqcolon ((U_R^n, V_R^n), (U_L^n, V_L^n)).
\end{equation}
We note that the variables $V_R^n=V_R^n(U_{-1}^n, U_{0}^n)$, $V_L^n=V_L^n(U_{-1}^n, U_{0}^n)$ depend only on the cell averages of the state variable $U$.

\begin{remark}
  As it holds
  \[T^i_j[U^n] - T^i_{j-1}[U^n] = \int_0^1 A_i(\Phi(s; U_{j-1}^n, U_j^n))\dds \Phi(s; U_{j-1}^n, U_j^n) \, ds \]
for both $i=1$ and $i=2$ the relaxed scheme does not require the evaluation of non-local operators away from the interface. Also the direction of integration within $T^i$ is not relevant when updating the numerical solution away from the interface. When imposing path conservative Kirchhoff conditions the coupling data for the auxiliary variable can be computed in terms of $\Sigma^{n,-}=V_R^n - T_{-1}^1[U^n]$ and $\Sigma^{n,+}=V_L^n - T_{0}^2[U^n]$, see Section~\ref{sec:kirchhoff}, which avoids the computation of the sums \eqref{eq:T1discrete} and \eqref{eq:T2discrete} in practice.
\end{remark}

\begin{proposition}\label{prop:pathconservativity}
  Let $\tilde \Phi$ denote the reverse family of paths to $\Phi$, i.e., $\widetilde \Phi(s; U_1, U_2)= \Phi(1-s; U_2, U_1)$ for $s\in[0,1]$. Suppose that $A_1=A_2\eqqcolon A$ and the vectors $U_a^1$ and $U_a^2$ are such that
  \begin{multline}\label{eq:consistencycondition}
    \int_0^1 A(\widetilde \Phi(s; U_1, U_a^1)) \dds \widetilde \Phi(s; U_1, U_a^1) \, ds + \int_0^1 A(\widetilde \Phi(s; U_a^2, U_2)) \dds \widetilde \Phi(s; U_a^2, U_2) \, ds\\
    = \int_0^1 A(\Phi(s; U_1, U_2)) \dds \Phi(s; U_1, U_2) \, ds
  \end{multline}
  for any $U_1\in \Omega_1$ and $U_2\in \Omega_2$. Suppose further that the coupling condition \eqref{eq:kirchhoff} is imposed to system \eqref{eq:couplednonconservative}. Then scheme \eqref{eq:relaxedcouplingscheme} is $\Phi$-conservative and identical to the uncoupled scheme~\eqref{eq:relaxedscheme}.
\end{proposition}
\begin{proof} It is sufficient to show that the scheme is $\Phi$-conservative at the interface. Making use of \eqref{eq:modVR} and \eqref{eq:modVL} we find that $D^{n,+}_{-1/2} = \frac 1 2 (\widetilde V_R - P_1(U_0^-)) + \frac 12 \Sigma^-$ and $D^{n,-}_{-1/2}=-\frac 1 2 (\widetilde V_L + P_2(U_0^+)) - \frac 1 2 \Sigma^+$ and therefore using \eqref{eq:rssystem1} and \eqref{eq:rssystem2} we obtain
  \begin{align*}
    D^{n,+}_{-1/2} + D^{n,-}_{-1/2} &= \frac 1 2 (\widetilde V_R - \widetilde V_L) - \frac 1 2 (P_1(U_0^-) + P_2(U_0^+)) + \frac 1 2 (\Sigma^- - \Sigma^+) \\
    &= - (P_1(U_0^-) + P_2(U_0^+)).
  \end{align*}
  As the last expression coincides with the left-hand side in \eqref{eq:consistencycondition} for $U_1=U_0^-$ and $U_2=U_0^+$ path consistency at the interface follows. 
\end{proof}
We note that Proposition~\ref{prop:pathconservativity} also holds in case of coupling conditions of the form \eqref{eq:modkirchhoff}.

\begin{remark}\label{rem:bc}
  In practice scheme~\eqref{eq:relaxedcouplingscheme} is used to compute numerical solutions on two bounded domains connected at an interface and problem specific boundary data is required. For this purpose suitable ghost data can be employed. Assuming a discretization by \eqref{eq:relaxedcouplingscheme} for $j=-N,\dots,N-1$ this is realized by specifying the data $U_{-(N+1)}^n$ and $U_N^n$. In consistency with the continuous zero-relaxation limit we assume then $U_j^n=U_a^1$ for $j<-N-1$ and $U_j^n=U_a^2$ for $j>N$, respectively. Thus the discretized auxiliary variable is given as $V_j^n = 0$ for both $j<-N-1$ and $j>N$ as well as
  \begin{align*}
    V_{-(N+1)}^n(U^n_{-(N+1)}) = T^1_{-(N+1)}[U^n] &= \int_0^1 A_1(\Phi(s; U_a^1, U_{-(N+1)}^n)) \dds \Phi(s; U_a^1, U^n_{-(N+1)}) \, ds, \\
    V_{N}^n(U_{N}^n) = T^2_{N}[U^n] &= -\int_0^1 A_2(\Phi(s; U_N^n, U_{a}^2)) \dds \Phi(s; U_N^n, U_{a}^2) \, ds.
  \end{align*}
\end{remark}

To discuss conservation properties we proceed in analogy to the derivation in Section~\ref{sec:relaxedscheme} and write scheme \eqref{eq:relaxedcouplingscheme} in the form
  \begin{equation}\label{eq:coupledconservativeform}
  U_j^{n+1} = U_j^n - \frac{\Delta t }{\Delta x} \left( H_{j+1/2}^{n,-} -  H_{j-1/2}^{n,+}\right) \quad \text{for all }j \in \Z.
\end{equation}
Away from the interface the numerical fluxes are then given by
\begin{align*}
   H_{j-1/2}^{n,+} = H_{j-1/2}^{n,-} =  \frac 1 2 \, (T^i_{j-1}[U^n] + T^i_j[U^n])  - \frac 1 2  \sqrt{\Lambda_i} (U_j^n-U_{j-1}^n)
\end{align*}
for $j\in \Z\setminus \{0\}$ and $i$ chosen according to the sign of $j$. Given the coupling data \eqref{eq:RSlimit} the numerical fluxes at the interface take the form
\begin{align*}
  H_{-1/2}^{n,-} &=  \frac 1 2 \, (T_{-1}^1[U^n] + V_{R}^n)  - \frac 1 2  \sqrt{\Lambda_1} (U_R^n-U_{-1}^n), \\
  H_{-1/2}^{n,+} &=  \frac 1 2 \, (V_{L}^n + T_0^2[U^n])  - \frac 1 2  \sqrt{\Lambda_2} (U_0^n-U_{L}^n).
\end{align*}
\begin{proposition}
  Suppose that $m_1=m_2$ and that the systems in \eqref{eq:couplednonconservative} on both axes are conservative so that $A_1=DF_1(U)$ and $A_2=DF_2(U)$. As in Remark~\ref{rem:bc} we consider two discretized and coupled bounded domains.
  \begin{enumerate}[label=(\alph*)]
  \item Suppose the coupling condition~\eqref{eq:kirchhoff} is imposed. If additionally the flux functions satisfy $F_1(U_a^1)=F_2(U_a^2)$ then scheme~\eqref{eq:relaxedcouplingscheme} is conservative at the interface, i.e., $ H_{-1/2}^{n,-} =  H_{-1/2}^{n,+}$.
  \item Under the coupling condition \eqref{eq:modkirchhoff} the scheme is conservative iff the boundary fluxes satisfy
    \[
      H_{N-1/2}^{n, -} - H_{-N-1/2}^{n, +} = F_1(U_a^1) - F_2(U_a^2).
      \]
  \end{enumerate}
\end{proposition}
 \begin{proof}
Due to~\eqref{eq:laxcondition} it holds $H_{-1/2}^{n,-} = V_R^n$ and $H_{-1/2}^{n,+}= V_L^n$, see also~\cite{herty2023centr}. As in the conservative case the coupling condition \eqref{eq:kirchhoff} implies $V_R + F_1(U_a^1) = V_L + F_2(U_a^2)$ statement (a) follows from the assumptions on the flux functions. 

Let us denote the total mass of the numerical solution at time $t^n$ by $M^n\coloneqq \sum_{j=-N}^{N-1} U_j^n$. From~\eqref{eq:coupledconservativeform} we obtain for the difference concerning two consecutive time instances
\[
  M^{n+1}-M^n = \frac{\Delta t}{\Delta x}\left( H_{-N-1/2}^{n,+} - H_{N-1/2}^{n, -} \right) - \frac{\Delta t}{\Delta x}\left( H_{-1/2}^{n,-} - H_{-1/2}^{n, +} \right).
\]
Noting then that under the coupling condition \eqref{eq:modkirchhoff} it holds $H_{-1/2}^{n,-} - H_{-1/2}^{n, +} = V_L^n - V_R^n = F_1(U_a^1) - F_2(U_a^2)$ statement (b) follows.
\end{proof}

\section{Numerical experiments}\label{sec:num}
In this section we present numerical experiments for the relaxed scheme in an uncoupled and a coupled setting. The uncoupled experiment is conducted over a bounded interval discretized over $N$ cells with imposed homogeneous Neumann boundary conditions. In the coupled case we discretize two bounded coupled intervals as described in Remark~\ref{rem:bc}. The parameter matrix $\Lambda$ is chosen as $\mu I$ so that $\mu$ provides a suitable bound over the squared eigenvalues, see Proposition~\ref{prop:mu}. We use segment paths in the numerical computation, i.e., the family $\Phi$ is such that $\Phi(s; W_1, W_2) = W_1 + s (W_2 - W_1)$ for all $W_1, W_2 \in \Omega$ and $s\in[0,1]$. The time increment is computed by the formula
\begin{equation}\label{eq:cfl}
\Delta t = \text{CFL} \frac{\Delta x}{ \sqrt{\mu}}.
\end{equation}
Details on the number of mesh cells and the Courant number used are provided in the individual experiment descriptions. The employed computer programs are implemented in the Julia programming language~\cite{bezanson2017julia}.

\subsection{Two-layer shallow water system}\label{sec:5.1}
To study the behavior of the uncoupled relaxation scheme~\eqref{eq:relschemeconservative} in the relaxation limit we use the two-layer shallow water system as an example of a nonconservative system. The model governs the flow of two superimposed layers of fluids over a flat bottom topography in a one-dimensional channel~\cite{castro2001}.
It can be stated in the form \eqref{eq:nonconservativesystem} with state vector and system matrix given by
\begin{equation}
  U = \begin{pmatrix}
        h_1\\
        q_1\\
        h_2\\
        q_2
      \end{pmatrix}, \qquad
      A(U) = \begin{pmatrix}
            0 & 1 & 0 & 0\\
            -\frac{q_1^2}{h_1^2} + g h_1 & \frac{2 q_1^2}{h_1^2}  & g h_1 & 0 \\
            0 & 0 & 0 & 1\\
            r g h_2 & 0 & -\frac{q_2^2}{h_2^2} + g h_2 & \frac{2 q_2^2}{h_2^2}
          \end{pmatrix}.
        \end{equation}
        The variables $h_i$ and $q_i$ denote the height and the mass flow of the upper and the lower layer represented by the index $i\in\{1,2\}$. Both layers are assumed to have constant density $\rho_i$ and $r=\rho_1/\rho_2$ represents their ratio, which we set to $0.9$. Furthermore, $g=9.81$ refers to the gravitational constant. We note that it is due to the entries $a_{23}(U)$ and $a_{14}(U)$ within $A(U)$, which relate the two layers, that the system is non-conservative. The spectral radius of the matrix $A(U)$ for our choice of $r$ is approximately given by
        \begin{equation}\label{eq:swspectral}
          \rho(A(U)) \approx \frac{|q_1 + q_2|}{h_1 + h_2} + \sqrt{g (h_1 + h_2)}.
        \end{equation}
\begin{figure}
  \centering
  \includegraphics{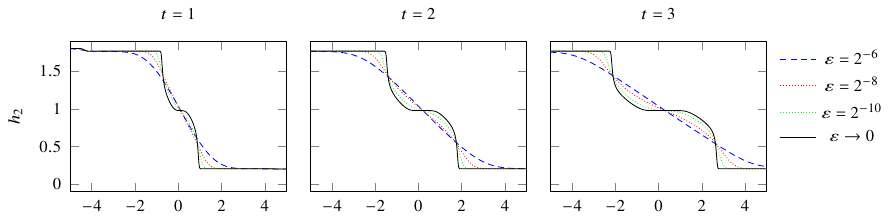}
    \caption{Height of the second layer $h_2$ in the internal dam break experiment over the computational domain in various computations using the relaxation (colored lines) and the relaxed scheme (black lines). The first layer is given by $h_1 = 2 - h_1$. As the relaxation rate decreases the numerical results of the relaxation scheme approach the ones of the relaxed scheme.}\label{fig:sw}
  \end{figure}
        
        We consider an experiment from \cite{castro2007well} modeling an internal dam break. The initial conditions in this case are
        \begin{equation}\label{eq:swini}
          q_1^0=q_2^0\equiv 0, \qquad h_1^0(x) = \begin{cases}
            0.2 &\text{if }x<0 \\
            1.8 &\text{if }x>0
          \end{cases},
          \qquad
          h_2^0(x) = \begin{cases}
            1.8 &\text{if }x<0 \\
            0.2 &\text{if }x>0
          \end{cases}
        \end{equation}
        and the considered spatial domain is $(-5, 5)$. In view of \eqref{eq:swspectral} and the initial data we choose $\mu=25$ in the CFL condition \eqref{eq:cfl}. The path integrals occurring in the schemes are approximated using the five point Gauss--Lobatto quadrature.
        Figure~\ref{fig:sw} shows numerical solutions for the relaxation scheme~\eqref{eq:relschemeconservative} and the limit scheme~\eqref{eq:relaxedscheme} on 4000 mesh cells in terms of height of the second layer $h_2$ for Courant number $\text{CFL}=0.9$. While for larger $\varepsilon$ the relaxation scheme predicts a straight interface between the two fluids the interface becomes more detailed as $\varepsilon$ decreases clearly approaching the prediction of the relaxed scheme as $\varepsilon \to 0$.

\begin{table}\scriptsize
  \centering
  \caption{Errors and EOCs with respect to the relaxation rate (left) and the number of mesh cells (right) for the two-layer shallow water system.}\label{tab:sw}
  \begin{tabular}{lllll}
    \toprule
$\varepsilon$ & $E_{\text{rel}, h_1}$ & EOC & $E_{\text{rel}, h_2}$ & EOC \\ \midrule
$2^{-7}$ & $1.66 \times 10^{-1}$ &  & $1.58 \times 10^{-1}$ &  \\
$2^{-8}$ & $8.80 \times 10^{-2}$ & 0.92 & $8.34 \times 10^{-2}$ & 0.92 \\
$2^{-9}$ & $4.52 \times 10^{-2}$ & 0.96 & $4.28 \times 10^{-2}$ & 0.96 \\
$2^{-10}$ & $2.28 \times 10^{-2}$ & 0.98 & $2.16 \times 10^{-2}$ & 0.99 \\
    $2^{-11}$ & $1.15 \times 10^{-2}$ & 0.99 & $1.08 \times 10^{-2}$ & 0.99 \\ \bottomrule
  \end{tabular} \hspace{1em}
    \begin{tabular}{rllll}
    \toprule
N & $E_{N, h_1}$ & EOC & $E_{N, h_2}$ & EOC \\ \midrule
500 & $4.55 \times 10^{-2}$ &  & $4.25 \times 10^{-2}$ &  \\
1000 & $2.39 \times 10^{-2}$ & 0.93 & $2.23 \times 10^{-2}$ & 0.93 \\
2000 & $1.21 \times 10^{-2}$ & 0.98 & $1.12 \times 10^{-2}$ & 0.99 \\
4000 & $6.15 \times 10^{-3}$ & 0.98 & $5.69 \times 10^{-3}$ & 0.98 \\
8000 & $3.07 \times 10^{-2}$ & 0.99 & $2.83 \times 10^{-3}$ & 1.00 \\ \bottomrule
  \end{tabular}
\end{table}
        
At the time instance $t=0.33$ we consider the relaxation error $E_{\text{rel}, h_i}$ that compares for various relaxation rates the numerical solution of the relaxation scheme with the one of the relaxed scheme in terms of the variable $h_i$ both over $4000$ mesh cells in the $L^1$ norm. The numerical solutions used in the error computations have been computed using the reduced Courant number $\text{CFL}=0.1$ and smooth initial data, in which $h^0_2(x) =  0.2 + \frac{1.6}{1 + \exp(-5x)}$ and $h_1^0(x) = 2-h_1^0(x)$ replace the discontinuous profile in~\eqref{eq:swini}. Table~\ref{tab:sw} shows a decrease of this error with $\varepsilon$ regarding both components $h_1$ and $h_2$. The experimental order of convergence (EOC\footnote{The EOC is computed by the formula $\text{EOC} = \log_2(E_1/E_2)$ with $E_1$ and $E_2$ denoting the error in two consecutive lines of the table.}) also provided in the table suggest first order convergence with respect to the relaxation rate. Using the numerical results we also present a grid convergence study with respect to the spatial error of the relaxed scheme and the variable $h_i$ that we denote by $E_{N,h_i}$ in Table~\ref{tab:sw}. Again the error is considered with respect to the $L^1$ norm at time instance $t=0.33$ and the results indicate first order convergence.

\subsection{Coupled blood flow model}\label{sec:blood}
In this section we apply our coupling approach to a model of blood flow through the human arterial system. While this process can be modeled in great detail using three dimensional fluid structure interaction models based on the Navier--Stokes equations reduced lumped parameter models on one-dimensional networks can accurately describe the interaction between pressure waves and the vessel geometry with reduced computational cost, see e.g.~\cite{MR2500548, quarteroni2004mathem} for details on the models. We consider such a reduced model studied in~\cite{formaggia2003one} given by the system
\begin{equation}\label{eq:bloodmodel}
  \begin{split}
  \ddt a + \ddx (au) &=0, \\
    \ddt u + (2 \alpha - 1) u \ddx u + (\alpha - 1) u^2 \ddx a + \rho^{-1} \ddx p_i &= - K_R \frac{u}{a},
  \end{split}
\end{equation}
where the state variables $a$ and $u$ represent the section area and the axial velocity, respectively. The parameter $K_R$ is used to model the viscosity of blood and is chosen as $8\pi \times 10^{-4}$ and $\rho$ refers to the blood density, which we set to $1$ for simplicity. We further assume that the pressure of the system is given by
\[
  p_i =  \beta_i \left( \sqrt{a} - \sqrt{a_0}\right), \qquad  \beta_i= \frac{E_i h_0 \sqrt{\pi}}{a_0} 
\]
with $a_0$ denoting the reference area of the vessel, $h_0$ the thickness of the vessel wall and $E_i$ the Young modulus corresponding to the elasticity of the vessel. We consider the blood flow through two connecting vessels for which the Young modulus and therefore the pressure function differ, which we realize by applying the coupling approach developed in Section~\ref{sec:coupling}.

Model~\eqref{eq:bloodmodel} can be written as a nonconservative system of the form \eqref{eq:nonconservativesystem} with an additional source term, so that state vector and system matrix are given by
\begin{equation}\label{eq:bloodmatrix}
  U = \begin{pmatrix}
        a \\ u
      \end{pmatrix}, \qquad
  A_i(U) = \begin{pmatrix}
             u & a\\
             (\alpha - 1) u^2 + \frac{\beta_i}{2 \rho \sqrt{a}}& (2\alpha - 1) u
           \end{pmatrix}.
         \end{equation}
         The Coriolis coefficient $\alpha$ plays an important role: when taking $\alpha=1$ as proposed in \cite{smith2002} system \eqref{eq:bloodmodel} becomes conservative and can be written as a balance law with flux function $F_i(U)=(au,~ \frac{1}{2} (2\alpha-1)u^2 + \rho^{-1}p_i)^T$. In fact, coupling conditions that are used to model connecting, branching and merging vessles in model~\eqref{eq:bloodmodel} are mostly motivated from the conservative case $\alpha=1$, cf.~\cite{formaggia2003one}.

           In our numerical experiment we couple model~\eqref{eq:bloodmodel} by considering \eqref{eq:bloodmatrix} in the coupled setting \eqref{eq:couplednonconservative} (additionally considering the corresponding source terms). At the interface we impose the path-conservative Kirchhoff condition~\eqref{eq:kirchhoff} and note that in the case $\alpha=1$ this choice coincides with the common coupling strategy in the literature, where equality of the flow rates and the total pressure is imposed at the interface~\cite{quarteroni2004mathem}. For the vessel properties we choose the parameters $a_0=5$, $h_0=0.05$ and assume that the right vessel is more elastic than the right vessel by taking $E_1=0.5$ and $E_2=0.1$. The experiment can be compared to blood flow from an artificial graft from a vascular bypass to a vein in the vasulature. While we impose homogeneous Neumann boundary conditions at the right boundary we simulate a heart by prescribing the time dependent boundary pressure profile $P_v(t) = P_0 \sin( \frac{\pi}{2} (t - \frac 12))$ for $P_0=2 \times 10^{-3}$ at the left boundary. Corresponding boundary data in terms of the state variables is then obtained by extrapolation of the outgoing Riemann invariant of the conservative system, see~\cite{formaggia2006numer} for details. As initial data we take $a_0 \equiv 5$ and $u_0 \equiv 0$. For the design of the RS at the interface we take $U_a^1=U_a^2=(5,0)^T$. To solve the nonlinear system given by~\eqref{eq:rssystem1} and \eqref{eq:rssystem2} Newton's method is used with starting point $\Sigma^-=\Sigma^+=0$. A visualization of the root functions in various relevant cases has indicated uniqueness of the numerically obtained root. We rely on the explicit computation of the path integral in Appendix~\ref{sec:bloodpaths} instead of using a quadrature formula in the computations.
         
\begin{figure}[t]
  \centering
  \includegraphics{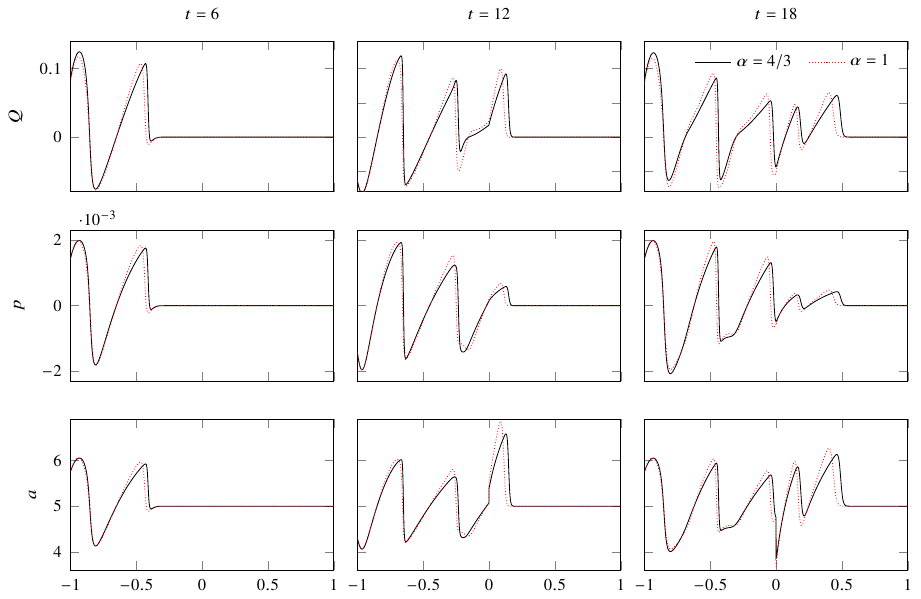}
    \caption{Numerical solutions for the coupled blood flow model in the nonconservative ($\alpha=4/3$) and the conservative ($\alpha=1$) case showing flow rate $Q$, pressure $p$ and section area $a$ at three time instances over the computational domain with coupling interface at $x=0$, where two blood vessels with different elasticity connect.}\label{fig:blood}
  \end{figure}  

We present numerical solutions on $4000$ mesh cells in Figure~\ref{fig:blood}. The computations have been obtained by the relaxed scheme~\eqref{eq:relaxedcouplingscheme} using the parameter $\mu=0.16$ and Courant number $\text{CFL}=0.9$ in~\eqref{eq:cfl}. Results in terms of flow rate $Q=au$, pressure and section area are shown under the assumption of Poiseuille flow ($\alpha=4/3$) and in the conservative case ($\alpha=1$). As the pulses propagate from the first to the second vessel the variation with respect to the pressure significantly decreases while variations from the reference section area increase, which also results in an increase of the flow rate. Notably, both flow rate and pressure stay continuous at the interface, whereas a jump in the section area occurs. In case of Poiseuville flow the pressure waves exhibit slightly lower magnitudes and slightly faster velocities when compared to the conservative case. We note that further numerical tests, in which we varied the family of paths and the truncation states $U_a^1$ and $U_a^2$, have not lead to qualitatively different numerical solutions. 

\begin{table}\footnotesize
  \centering
  \caption{Mesh convergence of the coupling errors with respect to \eqref{eq:continuouscoupling} and the $L^1$ errors for the coupled blood flow model~\eqref{eq:bloodmodel}. $N$ denotes the number of mesh cells per coupled domain.}\label{tab:blood}
  \begin{tabular}{rlrlrlrlr}
    \toprule
    $N$ & $E_{\Psi, 1}$ & EOC & $E_{\Psi, 2}$ & EOC & $E_{N,1}$ & EOC & $E_{N,2}$ & EOC \\  \midrule
    500 & $1.46 \times 10^{-3}$ &  & $2.12 \times 10^{-5}$ &  & $1.37 \times 10^{-1}$ &  & $1.83e-03$ &  \\
1000 & $7.88 \times 10^{-4}$ & 0.89 & $1.11 \times 10^{-5}$ & 0.93 & $8.37 \times 10^{-2}$ & 0.71 & $1.17 \times 10^{-3}$ & 0.65 \\
2000 & $4.09 \times 10^{-4}$ & 0.95 & $5.59 \times 10^{-6}$ & 0.99 & $4.58 \times 10^{-2}$ & 0.87 & $6.79 \times 10^{-4}$ & 0.78 \\
4000 & $2.09 \times 10^{-4}$ & 0.97 & $2.70 \times 10^{-6}$ & 1.05 & $2.39 \times 10^{-2}$ & 0.94 & $3.73 \times 10^{-4}$ & 0.87 \\
8000 & $1.05 \times 10^{-4}$ & 0.98 & $1.22 \times 10^{-6}$ & 1.15 & $1.22 \times 10^{-2}$ & 0.97 & $1.97 \times 10^{-4}$ & 0.92 \\ \bottomrule
  \end{tabular}
\end{table}

Since in the continuous case the terms $T^{1, 0^-}[U](0^-)$ and $T^{2, 0^+}[U](0^+)$ are each given by a single path integral, coupling condition \eqref{eq:modkirchhoff} implies
\begin{equation}\label{eq:continuouscoupling}
  \begin{split}
  \int_0^1 &A_1(\Phi(s; U_a^1, U(0^-))) \dds \Phi(s; U_a^1, U(0^-)) \, ds  \\
    &\quad +\int_0^1 A_2(\Phi(s; U(0^+), U_a^2)) \dds \Phi(s; U(0^+), U_a^2) \, ds  = 0
  \end{split}
\end{equation}
for a.~e. $t>0$. We consider the left-hand side of~\eqref{eq:continuouscoupling} and let $E_{\Psi, 1}$ and $E_{\Psi, 2}$ denote the absolute values of its two components at time instance $t=12$. In Table~\ref{tab:blood} we present the behavior of these coupling errors under mesh refinement in the above experiment for $\alpha=4/3$. To reduce the effect of temporal errors very small time increments (by imposing Courant number $\text{CFL}=0.02$) have been used in the computation of the numerical solutions. The coupling errors significantly decrease as the mesh is refined and the computed EOCs indicate a convergence linear in the mesh width. In the same table we present the $L^1$ errors of the full numerical solution with respect to both system components at the same time instance, which we denote by $E_{N,1}$ and $E_{N,2}$. Again, the EOCs indicate convergence of first order for these errors.

\section{Conclusion}
We have introduced a new relaxation approach for coupled nonconservative hyperbolic systems relying on the framework in~\cite{dal1995defin}. The system relaxes towards a nonlocal Borel-measure of the solution, which depends on a chosen path family. The asymptotic expansion in Theorem~\ref{thm:expansion} verifies the correct relaxation limit and motivates a stability condition in line with the subcharacteristic condition necessary for well-posedness of relaxation of conservative systems. Discretizing the relaxation system by an implicit-explicit asymptotic-preserving scheme, we have recovered the path-conservative Lax--Friedrich scheme in the relaxation limit. A numerical experiment considering the two-layer shallow water system indicates first order of convergence of the relaxation scheme to the limit scheme with respect to the relaxation rate.

The relaxation system has allowed us to generalize the approach from~\cite{herty2023centralsystems} to couple two nonconservative systems. The direction of integration indicating the support of the Borel measure plays a crucial role in the formulation of the coupling problem. We propose to integrate towards the coupling interface, which allows for a notion of (linear) Riemann solvers consistent with the conservative case. Consistency at the interface with the relaxation approach requires a coupling condition in terms of the nonlocal Borel-measure. To handle local coupling conditions we propose an approach based on truncation within the operators~\eqref{eq:T1} and \eqref{eq:T2}. In this way we introduce path-conservative Kirchhoff conditions, for which we provide  a fully discrete scheme. We present an application to a model of blood flow in the vasculature, where our approach leads to new coupling conditions for the nonconservative case. 

{\small \subsection*{Acknowledgments} 
  The authors thank the Deutsche Forschungsgemeinschaft (DFG, German Research Foundation) for the financial support under Germany’s Excellence Strategy EXC-2023 Internet of Production 390621612 and under the Excellence Strategy of the Federal Government and the Länder, 333849990/GRK2379 (IRTG Hierarchical and Hybrid Approaches in Modern Inverse Problems), 320021702/GRK2326, 442047500/SFB1481 within the projects B04, B05  and B06, through SPP 2410 Hyperbolic Balance Laws in Fluid Mechanics: Complexity, Scales, Randomness (CoScaRa) within the Project(s) HE5386/26-1 and HE5386/27-1, MU 1422/9-1 and through SPP 2298 Theoretical Foundations of Deep Learning within the Project(s) HE5386/23-1, Meanfield Theorie zur Analysis von Deep Learning Methoden (462234017). Support by the EU DATAHYKING No. 101072546 as well as by the ERS Open Seed Fund of RWTH Aachen University through project OPSF781 is also acknowledged.  
}

\bibliographystyle{abbrvurl} 
\bibliography{references.bib}

\appendix
\section{Path integrals for the blood flow model}\label{sec:bloodpaths}
Let $A_i$ be given by \eqref{eq:bloodmatrix} and $\Phi$ denote the family of segment paths. We take a generic path $\phi(s)\defeq \Phi(s; U_1, U_2)$ with components $\phi_1$ and $\phi_2$, where $U_1=(a_1, u_1)^T$ and $U_2=(a_2, u_2)^T$. Then the path integral \[\int_0^1 A_i(\phi(s))\dds \phi(s) \, ds  \]
 is given by
 \begin{align*}
   \int_0^1 \phi_2(s) \, ds \, \phi_1^\prime + \int_0^1 \phi_1(s) \, ds \, \phi_2^\prime
   &= \int_{u_1}^{u_2} u \, du \, \frac{\phi_1^\prime}{\phi_2^\prime} + \int_{a_1}^{a_2} a \, da \, \frac{\phi_2^\prime}{\phi_1^\prime} \\
   & = \frac 1 2 \left( (u_2 + u_1) (a_2 - a_1) + (a_2 + a_1) (u_2 - u_1)\right) 
 \end{align*}
 in the first component and by
 \begin{align*}
    2 (\alpha -1) &\int_0^1 \phi_2(s) \, ds \, \phi_1^\prime + \int_0^1 \frac{\beta_i}{2 \rho \sqrt{\phi_1(s)}} \, ds \, \phi_1^\prime + (2 \alpha - 1) \int_0^1 \phi_2(s) \, ds \, \phi_2^\prime \\
   &= 2 (\alpha -1) \int_{u_1}^{u_2} u^2 \, du \,  \frac{\phi_1^\prime}{\phi_2^\prime} + \int_{a_1}^{a_2} \frac{\beta_i}{2 \rho \sqrt{a}} \, da + (2 \alpha - 1) \int_{u_1}^{u_2} u \, du \\
   &= \frac{2(\alpha-1)}{3} (u_2^3 - u_1^3)\,  \frac{a_2 - a_1}{u_2 - u_1} + \frac{\beta_i}{\rho} (\sqrt{a_2} - \sqrt{a_1}) + \frac{2\alpha -1}{2}(u_2^2 - u_1^2)
 \end{align*}
 in the second component. 
\end{document}